\documentclass{amsart}
\usepackage{amsrefs}
\usepackage{cleveref}
\usepackage{upref}

\newtheorem{theorem}{Theorem}
\newtheorem{proposition}{Proposition}
\newtheorem{lemma}{Lemma}

\newtheorem{claim}{Claim}

\def \ue {u_\epsilon}
\def \le {\lambda_\epsilon}
\def \he {h_\epsilon}

\def \hfe {\hat{f}_\epsilon}
\def \hue {\hat{u}_\epsilon}
\def \hhe {\hat{h}_\epsilon}
\def \hge {\hat{g}_\epsilon}
\def \ye {y_\epsilon}
\def \theps {\tilde{h}_\epsilon}
\def \tfe {\tilde{f}_\epsilon}
\def \tue {\tilde{u}_\epsilon}

\def \we {w_\epsilon}
\def \xe {x_\epsilon}
\def \me {\mu_\epsilon}
\def \rr {\mathbb{R}}
\def \rn {\mathbb{R}^n}
\def \eps {\epsilon}
\def \critk {2^\star(k)}
\def \crit {2^\star}
\def \mbar {\overline{M}}
\title[Concentration analysis for nonlinear equations]{Concentration analysis for elliptic critical equations with no boundary control: ground-state blow-up}

\author{Hussein Mesmar}
\address{Hussein Mesmar, Institut \'Elie Cartan, Universit\'e de Lorraine, BP 70239, F-54506 Vand\oe uvre-les-Nancy, France}
\email{houssein.mesmar91@hotmail.com}

\author{Fr\'ed\'eric Robert}
\address{Fr\'ed\'eric Robert, Institut \'Elie Cartan, Universit\'e de Lorraine, BP 70239, F-54506 Vand\oe uvre-l\`es-Nancy, France}
\email{frederic.robert@univ-lorraine.fr}

\date{October 2nd, 2023. MSC: 35J60, 58J10}

\begin{document}
\begin{abstract} We perform the apriori analysis of solutions to critical nonlinear elliptic equations on manifolds with boundary. The solutions are of minimizing type. The originality is that we impose no condition on the boundary, which leads us to assume $L^2-$concentration. We also analyze the effect of a non-homogeneous nonlinearity that results in the fast convergence of the concentration point.
\end{abstract}
\maketitle
\centerline{\it Dedicated to Yihong Du on the occasion of his 60th birthday}

\section{Introduction}
\subsection{Context and main results} Let $(M,g)$ be a Riemannian manifold of dimension $n\geq 3$, with or without boundary $\partial M$. When $\partial M\neq\emptyset$, $M$ denotes the interior of the manifold and $\overline{M}$ denotes its closure, so that $\overline{M}=M\cup\partial M$: in particular, $M$ is open in $\overline{M}$. We let $a,f\in C^0(\mbar )$ be functions and we consider $u\in C^2(M)$ solution to
\begin{equation}\label{eq:init}
\Delta_g  u + au  =f u^{\crit-1}\, ;\, u>0    \hbox{ in  } M.
\end{equation}
where $\Delta_g:=-\hbox{div}_g(\nabla)$ is the Laplacian with minus sign convention and $\crit:=\frac{2n}{n-2}$ is critical for the Sobolev embeddings $H_1^2(M)\hookrightarrow L^{\crit}(M)$. Here, the Sobolev space $H_1^2(M)$ is the completion of $\{u\in C^\infty(M)/\, \Vert u\Vert_{H_1^2}<\infty\}$ for the norm $\Vert\cdot\Vert_{H_1^2}:=\Vert\nabla\cdot\Vert_2+\Vert\cdot\Vert_2$. In the case of a Euclidean smooth domain $\Omega\subset \rn$, then $a,f\in C^0(\bar{\Omega})$ and we consider $u\in C^2(\Omega)$ solution to
\begin{equation}\label{eq:init:euc}
\Delta  u + au  =f u^{\crit-1}\, ;\, u>0    \hbox{ in  } \Omega.
\end{equation}
where  $\Delta:=-\hbox{div}(\nabla)$ is the Euclidean Laplacian. Due to the critical exponent $\crit$, there might be families of solutions to \eqref{eq:init} that are not relatively compact in $C^2_{loc}(M)$. For instance, given $x_0\in\rn$ and $\mu>0$, define the {\it Bubble} as 
\begin{equation*}
x\mapsto U_{\mu,x_0}(x):=\left(\frac{\mu}{\mu^2+\frac{c|x-x_0|^2}{n(n-2)}}\right)^{\frac{n-2}{2}}.
\end{equation*}
Then for any domain $\Omega\subset \rn$, $U_{\mu,x_0}$ is a solution to \eqref{eq:init:euc} when $a\equiv 0$ and $f\equiv c$. Moreover, if $x_0\in \bar{\Omega}$, then $\sup_{\Omega}U_{\mu, x_0}\to +\infty$ as $\mu\to +\infty$. In the Riemannian context, for $x\in M$ and $\mu>0$, the {\it Bubble} writes
\begin{equation*}
x\mapsto U^{(M,g)}_{\mu,x_0}(x):=\left(\frac{\mu}{\mu^2+\frac{cd_g(x,x_0)^2}{n(n-2)}}\right)^{\frac{n-2}{2}}.
\end{equation*}

\smallskip\noindent Concerning terminology, we say that a family $(\ue)_\eps\in C^0(M)$ blows-up if $\lim_{\eps\to 0}\Vert \ue\Vert_{\infty}=+\infty$.

\smallskip\noindent When $\partial M=\emptyset$, the description of blowing-up families of \eqref{eq:init} with bounded $L^{\crit}-$norm has been performed by Druet-Hebey-Robert \cite{dhr}. The main result in \cite{dhr} is that blowing-up families are controled from above by $\left(\sum_iU_{\mu_{i,\eps},x_{i,\eps}}\right)_\eps$, for given families $(x_{i,\eps})_\eps\in M$ and $(\mu_{i,\eps})_\eps\to 0$. With this control, it is possible to give informations on the localization of the limits $x_{i,\infty}:=\lim_{\eps\to 0}x_{i,\eps}$: see Druet \cite{druet:jdg}. This analysis extends to manifolds with boundary provided a boundary condition like Dirichlet (see Ghoussoub-Mazumdar-Robert \cite{gmr}) or Neumann (see Druet-Robert-Wei \cite{drw}). See Premoselli \cite{premoselli} for a more recent point of view.

\medskip\noindent The first objective of the present work is to perform an analysis similar to \cite{dhr} and \cite{druet:jdg} without condition on $\partial M\neq \emptyset$. Tackling such generality requires additional assumption: the relevant notion here is $L^2-$concentration that already appeared in Djadli-Druet \cite{dd} (see \eqref{hyp:conctration:L2:article:blowup} below).

\smallskip\noindent Our second objective is to analyse the effect of a nonconstant function $f$ in \eqref{eq:init}. In the case of a single peak, concentration occurs at a critical point. We prove that when this critical point is nondegenerate, then the family of concentration points converges very fast to its limit (see \eqref{d:est:petit:o:de:mu:eps} below): this does not generally happen for a constant function $f$. A similar control appears in Malchiodi-Mayer \cite{mm}.

\smallskip\noindent As was shown by Aubin \cite{aubin0problemes}, below a threshold, blow-up cannot occur. In this manuscript, we are considering solutions $(\ue)$ that carry the minimal energy for blow-up, namely ground-state type solution. The minimal energy is given by the best constant in Sobolev embeddings:
\begin{align} \label{ineq:sobolev}
    \frac{1}{K_0 (n)} = \inf_{\varphi \in D_1^2(\rr^n) \setminus \{0 \}} \frac{\int_{\rr^n} |{\nabla} \varphi |^ 2 \mathrm{d} X}{\left (\int_{\rr^n} | \varphi |^{\crit} \mathrm{d} X \right)^{\frac{2}{\crit}}},
\end{align}
where $D_1^2(\rr^n)$ is the completion of $C_c^{\infty} (\rr^n)$ for the norm $\varphi\mapsto \Vert\nabla\varphi\Vert_2$. Aubin \cite{aubin0problemes} and Talenti \cite{talenti0best} have computed this best constant and have showed that the extremals are exactly $C\cdot U_{\mu,x_0}$ for $C\neq 0$, $\mu>0$ and $x_0\in\rn$. Our main theorem for ground-state solutions is the following:
\begin{theorem}\label{theorem:Mesmar:1:article:blow:up}
Let $(M,g)$  be a smooth compact Riemannian manifold of dimension $ n\geq 4 $ with nonempty boundary $\partial M\neq\emptyset$. We fix $f\in C^2(\mbar )$ such that $f>0$. We consider a family $(\he)_\eps\in C^1(\mbar )$ and $f \in C^2(\mbar)$, $f>0$, such that there exists $h \in C^1(\mbar )$ and such that $\Delta_g+h $ is coercive and 
\begin{equation}\label{lim:h:f}
\lim_{\eps\to 0}\he=h \hbox{ in }C^1(\mbar ).
\end{equation} 
We let  $(\ue)_\eps\in C^2(\mbar ) $ be a family of solutions to  
\begin{equation}\label{principal:eq:critical:variete:with:boundary}
    \Delta_g \ue +  \he \ue= f \ue^{\crit - 1} \hbox{ in }M.
\end{equation}
Let $\xe\in \mbar $ and $\me>0$ be such that
 \begin{equation*}
     \ue ( \xe )   = \sup_{\mbar }  \ue  =  \me^{1-\frac{n}{2}} 
     \end{equation*}
We assume that
\begin{itemize}
\item $\ue\to  0$  in $L^2(M)$,
\item $\lim_{\eps\to 0}\xe=x_0\in M$ is an interior point of $M$,
\item The solution has minimal-type energy, that is
\begin{equation*}
\lim_{\eps\to 0}\int_M f \ue^{\crit}\, dv_g=\frac{1}{K_0(n)^{\frac{n}{2}} f (x_0)^{\frac{n-2}{2}}}
\end{equation*}
\item The Hessian $\nabla^2f(x_0)$ is nondegenerate.
\end{itemize}
Then $x_0$ is a critical point of $f$ and 
\begin{align}\label{d:est:petit:o:de:mu:eps}
    d_g (\xe, x_0) = o(\me)\hbox{ as }\eps\to 0,
\end{align}
and for all $\omega\subset M$ such that $\overline{\omega}\subset M$ and $\delta_0>0$, there exists $C(\omega,\delta_0)>0$ such that
\begin{equation*}
\ue(x)\leq C(\omega,\delta_0)\left(\frac{\me}{\me^2+d_g(x,x_0)^2}\right)^{\frac{n-2}{2}}+C(\omega,\delta_0)\sup_{\partial B_{\delta_0}(x_0)}\ue
\end{equation*}
for all $x\in \omega$. In addition, assuming that for all $ \delta> 0 $, we have that 
\begin{equation}\label{hyp:conctration:L2:article:blowup} 
    \lim_{\epsilon \to 0} \frac{\int_{M \setminus B (x_0 , \delta)}\ue^ 2\, dv_g }{\int_M u_ \epsilon^ 2 \, dv_g} = 0 \, \hbox{ for } \,  n\in \{4, 5, 6\},
 \end{equation}
then  
\begin{align}\label{relation:de:h0:avec:courbure:scalaire:cas:variete:a:bord}
    h (x_0) = \frac{n-2}{4(n-1)}\left(  \hbox{Scal}_g(x_0)- \frac{ n-4}{2}\cdot \frac{ \Delta_g f  (x_0)}{ f (x_0) } \right),
\end{align}
where $ \hbox{Scal}_g$ is the scalar curvature of $(M,g)$.
\end{theorem}
\noindent{\bf Remark: }{\it Theorem \ref{theorem:Mesmar:1:article:blow:up} applies to the case of a bounded domain of $\rn$ endowed with the Euclidean metric $g:=\hbox{Eucl}$. In this situation, $M=\Omega\subset\rn$ is a domain, $\Delta_g=-\sum_i\partial_{ii}$, $d_g(x,y)=|x-y|$ is the usual Euclidean norm for $x,y\in\rn$ and $\hbox{Scal}_g=0$.}

\smallskip\noindent The control \eqref{d:est:petit:o:de:mu:eps} is remarkable since it does not hold when $f$ is degenerate. Indeed, when $f\equiv 1$ there is an abundance of blowup profiles with various speeds of convergence of the $(\xe)$'s to their limit, see for instance Premoselli \cite{premo:jgea}.

\medskip\noindent The restriction of dimension $n\geq 4$ is not surprising: indeed, see Corollary 6.4 in Druet-Hebey \cite{AB}, \eqref{hyp:conctration:L2:article:blowup} does not hold in general for $n=3$. It is known since Aubin and Schoen that for $n=3$, blowup cannot be characterized by local arguments and involves global arguments, like the mass. In the general local context of Theorems \ref{theorem:Mesmar:1:article:blow:up}, no information is known regarding the boundary, which forbids to get any global information.

\subsection{Application to supercritical problems with symmetries} A natural set application of Theorem \ref{theorem:Mesmar:1:article:blow:up} is in the context of manifolds invariant under a group of isometries. We consider a compact Riemannian manifold $(X,g)$ of dimension $n\geq 3$, but without boundary ($\partial X=\emptyset$). The critical exponent $\crit$ can be improved by imposing invariance under the action of an isometry group. Let $G$ be a compact subgroup of isometries of $(X,g)$: we say that a function $u:X\to\rr$ is $G-$invariant if $u\circ \sigma=u$ for all $\sigma\in G$. It follows from Hebey-Vaugon \cite{hebey0sobolev} that the critical exponent in this setting is $\critk:=\frac{2(n-k)}{n-k-2}$ where $k:=\min_{x\in X}\hbox{dim }Gx$ and assuming that $1\leq k<n-2$. We refer to Hebey-Vaugon \cite{hebey0sobolev}, Saintier \cite{saintier0blow} and Faget \cite{faget2002best} for extensive considerations on problems invariant under isometries. In general, the quotient $X/G$ is not a manifold of dimension $n-k$. Following Saintier \cite{saintier0blow}, we make the following assumption on $G$:

\smallskip\noindent{\bf Assumption (H): }{\it For any $x_0\in X$ such that the orbit $Gx_0 $ is of dimension $k= \min\limits_{x \in X} dim Gx \geq 1$ and of volume $V_m=\min_{x\in X}\{\hbox{Vol}_g(Gx)/\, \hbox{dim }Gx=k\}$, there exists $ \delta> 0 $, and $ G '$ a closed subgroup of $ Isom_g (X) $ such that:
  \begin{enumerate}
  \item $ G'x_0 = Gx_0 $;
\item For all $ x \in B_{\delta} (Gx_0):= \lbrace y \in X / d_g (y; Gx_0) <\delta \rbrace $, then $ G'x $ is principal and $ G ' x \subset Gx $.
  \end{enumerate}
In particular $ dim\,G'x = dim\, Gx_0 = k $, $ \forall x \in B_{\delta} (Gx_0) $.}\par

\smallskip\noindent This assumption ensures that $ B_{\delta} (Gx_0) / G '$ is a Riemannian manifold of dimension $m:= n-k $ with a nontrivial boundary. In the sequel, for any $p\in\mathbb{N}$, we define   $C^p_G(X)$ as the space of $G-$invariant functions of $C^p(X)$.  We prove the following in the spirit of Faget \cite{faget:ade}.

\begin{theorem}\label{th:mesmar:2:supercitical:case:blow:up} Let $ (X, g) $ be a compact Riemannian manifold of dimension $ n $ without boundary, and let $G$ be a compact subgroup of isometries of $X$ which satisfies Assumption $(H)$ and such that $1\leq k<n-2$. Let $(\he)_\eps\in C^1_G(X)$ and $h \in C^1_G(X)$ be such that $\Delta_g+h $ is coercive and 
\begin{equation} \label{heps:intro:eng}
\lim_{\eps\to 0}\he=h >0\hbox{ in }C^1_G(X).
\end{equation}
Let $ (\ue)_\eps\in C^2_G(X) $ be a family of solutions to
\begin{equation} \label{principale:equation:supercritical:blowup}
    \Delta_g\ue + \he \ue =\le\ue^{\critk - 1}\, ;\, \ue>0 \hbox{ in }X, \;  \int_X\ue^{\critk}\, dv_g = 1
\end{equation}
We assume that
\begin{itemize}
\item $\ue\to 0$ strongly in $L^2(X)$,
\item The energy is of minimal type, that is
\begin{equation} \label{align:hyp:hme:blowup:supercritical}
 \lim_{\epsilon \to 0}\le = \frac{V_m^{1- \frac{2}{\critk} }}{K_0 (n-k)}, \hbox{ where }V_m:=\min_{x\in X}\{\hbox{Vol}_g(Gx)/\, \hbox{dim }Gx=k\}.
\end{equation}
\item For all point $z_0\in X$ such that $ \hbox{dim }Gz_0=k$ and $\hbox{Vol}_g(Gz_0)=V_m$, then the function
\begin{equation*}
\left\{\begin{array}{cccc}
\bar{v}: & B_\delta(Gz_0)/G' &\to & \rr\\
& G'x &\to & \hbox{Vol}_g(G'x)
\end{array}\right\}\hbox{ is nondegenerate at }Gz_0.
\end{equation*}
\end{itemize}
This latest assumption makes sense due to Assumption (H). Let $(\xe)_\eps\in X$ be such that $\ue(\xe)=\max_X\ue$ and define $\me^{-\frac{n-k-2}{2}}=\ue(\xe)$. Then there exists $x_0\in X$ such that  $\hbox{dim }Gx_0=k$ and $\hbox{Vol}_g(Gx_0)=V_m$ such that $\lim_{\eps\to 0}\xe=x_0$ and 
\begin{align} \label{distance:with:orbit:smal:o:of:mu}
      d_g (\xe, Gx_0) = o(\me).
\end{align}
Moreover, there exists $C>0$ such that
\begin{equation}\label{ineq:co:G}
\ue(x)\leq C\left(\frac{\me}{\me^2+d_g(x, Gx_0)^2}\right)^{\frac{n-2}{2}}+\left\{\begin{array}{cc}
o(\me)&\hbox{ if }n-k\geq 5\\
o\left(\me\sqrt{\ln\frac{1}{\me}}\right)&\hbox{ if }n-k=4
\end{array}\right.\end{equation}  
and 
\begin{equation} \label{id:scal:intro:eng}
h  (x_0) = \frac{n-k-2}{4 (n-k-1)} \left (\hbox{Scal}_{\bar{g}} (\bar{x}_0) +3 \frac{\Delta_{\bar{g}} \bar{v} (\bar{x}_0)}{\bar{v} (\bar{x}_0)} \right) \hbox{ when } \,  n-k \geq 4,
\end{equation}
where $\bar{g}$ is the metric on $B_\delta(Gx_0)/G'$ such that the canonical projection $(B_\delta(Gx_0),g) \to (B_\delta(Gx_0)/G',\bar{g})$ is a Riemannian submersion.
\end{theorem}

\section{Pointwise control}
We consider $(\ue)\in C^2(\mbar )$, $(\he)\in C^1(\mbar )$, $h\in C^1(\mbar)$, $f\in C^2(\mbar )$, $(\xe)_\eps\in M$ and $(\me)_\eps\in (0,+\infty)$ as in the statement of Theorem \ref{theorem:Mesmar:1:article:blow:up}. In the sequel, we let $i_g(M,x)>0$ be the injectivity radius of $(M,g)$ at an interior point $x\in M$.  
\begin{claim}\label{claim:cv:we} Set $\delta\in (0, i_g(M,x_0))$ and define
\begin{equation*}
     \we  (X) :=\me^{\frac{n-2}{2}} \ue (\exp_{\xe}(\me  X )) \hbox{ for any } X \in {B_{\frac{\delta}{\me }}(0)} \subset \rn
 \end{equation*}
Then 
\begin{equation}\label{cv:we:w}
\lim_{\epsilon \to 0}\we(X)=w(X)=\left(\frac{1}{1+\frac{f(x_0)|X|^2}{n(n-2)}}\right)^{\frac{n-2}{2}}\hbox{ for all }X\in \rn.
\end{equation}
Moreover, the convergence holds in  $C^2_{loc}(\mathbb{R}^n)$. In addition,
\begin{equation*}
           \lim\limits_{R\to + \infty } \lim\limits_{{\eps } \to 0} \int_{B_{\xe }(R \me)} f\ue^{\crit} \, dv_g  =\frac{1}{K_0(n)^{\frac{n}{2}} f (x_0)^{\frac{n-2}{2}}}.
       \end{equation*}
       In particular
\begin{equation}\label{estim:out:ball:bis}
           \lim\limits_{R\to + \infty } \lim\limits_{{\eps } \to 0} \int_{M\setminus B_{\xe }(R \me)} f\ue^{\crit} \, dv_g  =0.
       \end{equation}
       
\end{claim}
\noindent{\bf Proof of Claim \ref{claim:cv:we}:} We define the metric $g_\eps:=  \exp_{\xe}^\star g(\me\cdot)$ in $B_{\frac{\delta}{\me }}(0)\subset \rn$. Since, $ \me  \rightarrow 0 $ when $ {\mathbb {\epsilon}} \rightarrow 0 $, then $ g_\eps\to \xi$ in $C^2_{loc}(\rn)$ as $\eps\to 0$ where $\xi$ is the Euclidean metric. The function $\we$ satisfies the equation
    \begin{equation} \label{eq:twe} 
     \Delta_{g_\eps}  \we + \me^2 \tilde{h}_\epsilon \we  = \tfe \we ^{\crit-1}\hbox{ in }B_{\frac {\delta} \me } (0)
     \end{equation}
 where $ \theps(X)   = \he\left(\exp_{\xe  }(\me X)\right)$ and $\tfe (X)= f(\exp_{\xe}(\me  X ))$ for all $X\in B_{\frac {\delta} \me } (0)$. Since $0<\we\leq \we(0)=1$,  there exists $ w \in C^2 \left (\rr^n\right) $ such that the sequence $\we\to w$ in $C^2_{loc}(\rn)$ as $\eps\to 0$ up to extraction. Passing to the limit in \eqref{eq:twe}, we get that
 \begin{equation} \label{equ:w:critical:case}
      \Delta_\xi {w } = f(x_0)  {w }^{\crit-1}\hbox{ in }\rn,\, 0\leq w(0)=1.
      \end{equation}
It follows from Cafarelli-Gidas-Spruck \cite{cgs} that  $w(X)=\left( 1+\frac{f(x_0)|X|^2}{n(n-2)}\right)^{-\frac{n-2}{2}}$ for all $X\in \rn$. The change of variable $x=\hbox{exp}_{\xe}(\me X)$ yields
 \begin{align*}
    \int_{B_{\xe }(R \me)} f \ue^{\crit} \, dv_g  = \int_{B_R(0)} f(\hbox{exp}_{\xe }(\me  X) w_{\epsilon}^{\crit}  \mathrm {d} v_{ g_\eps }.
 \end{align*}
 Therefore, 
 \begin{eqnarray*}
 \lim_{R\to +\infty}\lim_{\eps\to 0}\int_{B_{\xe }(R \me)} f \ue^{\crit} \, dv_g & =& \lim_{R\to +\infty}\lim_{\eps\to 0}\int_{B_R(0)} f(\hbox{exp}_{\xe }(\me  x) w_{\epsilon}^{\crit}  \mathrm {d} v_{ g_\eps }\\
 &=&f(x_0)\int_{\rn}w^{\crit}\, dx=\frac{1}{K_0(n)^{\frac{n}{2}}f(x_0)^{\frac{n-2}{2}}},
 \end{eqnarray*}
 where we have used that $w$ is a solution to \eqref{equ:w:critical:case} and is an extremal for the Sobolev inequality \eqref{ineq:sobolev}. This proves Claim \ref{claim:cv:we}.\qed

 \begin{claim}\label{claim:2}  $\ue\to 0$ in $C^0_{loc}(M\setminus\{x_0\})$.
 \end{claim}
\noindent{\it Proof of the claim:} It follows from  \eqref{estim:out:ball:bis} and $f>0$ that for all $\delta>0$, we have that
\begin{equation*}
\lim_{\eps\to 0}\int_{M\setminus B_{x_0}(\delta)}\ue^{\crit}\, \, dv_g=0.
\end{equation*}
Let us fix  $\omega\subset M$ such that $\overline{\omega}\subset M\setminus\{x_0\}$. We let $\omega'$ open such that $\bar{\omega}\subset \omega'$ and $\bar{\omega'}\subset M-\{x_0\}$. Let $\eta\in C^\infty_c(\omega')$ such that $\eta(x)=1$ for all $x\in\omega$. Let us take $l>1$ to be fixed later. Integrating by parts as in Druet-Hebey (\cite{AB}, Theorem 6.1), we get that
\begin{eqnarray*}
\int_M\eta^2\ue^l\Delta_g\ue\, dv_g &=& \int_M\nabla(\eta^2\ue^l)\nabla\ue\, dv_g=\int_M l\eta^2\ue^{l-1}|\nabla\ue|^2\, dv_g+\int_M\nabla \eta^2\nabla\frac{\ue^{l+1}}{l+1} \, dv_g\\
&=&\frac{4l}{(l+1)^2}\int_M\eta^2|\nabla\ue^{\frac{l+1}{2}}|_g^2\, dv_g+\int_M\frac{\Delta\eta^2}{l+1}\ue^{l+1}\, dv_g
\end{eqnarray*}
Independently, for any $v\in C^1(M)$, integrating also by parts, we get that
$$\int_M(|\nabla (\eta v)|_g^2-\eta^2|\nabla v|_g^2)\, dv_g=-\int_M\eta v^2\Delta_g\eta\, dv_g.$$
Plugging these integrals together yields
\begin{eqnarray*}
&& \int_M \left| {\nabla}  \right( \eta \ue^{\frac{l+1}{2}}    \left)  \right|_g^2 \, dv_g  = \frac{(l+1)^2}{4l} \int_M \eta^2 \ue^l \Delta_g \ue \, dv_g  + \frac{l+1}{2l}\int_M \left( |{\nabla}  \eta |_g^2 + \frac{l-1}{l+1} \eta \Delta_g \eta \right) \ue^{l+1} \, dv_g  
 \end{eqnarray*}
We then get that
\begin{eqnarray*}
    && \int_M \left| {\nabla}  \right( \eta \ue^{\frac{l+1}{2}}    \left)  \right|_g^2 \, dv_g \\
    && = \frac{(l+1)^2}{4l} \int_M \eta^2 \ue^l \Delta_g \ue \, dv_g  + \frac{l+1}{2l}\int_M \left( |{\nabla}  \eta |_g^2 + \frac{l-1}{l+1} \eta \Delta_g \eta \right) \ue^{l+1} \, dv_g  \\
        &&=     \frac{(l+1)^2}{4l} \le   \int_M  \eta^2 \ue^{l+\crit-1} \, dv_g  - \frac{(l+1)^2}{4l} \int_M \eta^2 {h_{\eps }} \ue^{l+1} \, dv_g   \\
     & &+ \frac{l+1}{2l} \int_M \left( |{\nabla}  \eta |_g^2 + \frac{l-1}{l+1} \eta \Delta_g \eta \right) \ue^{l+1} \, dv_g   \\
      &&\leq   \frac{(l+1)^2}{4l} \le  \left ( \int_M \left( \eta \ue^{\frac{l+1}{2}} \right)^{\crit} \, dv_g  \right)^{\frac{2}{\crit}} \left( \int_{B_{x}(2\delta)} \eta \ue^{\crit} \, dv_g   \right)^{1-\frac{2}{\crit}}+ C\int_\Omega \ue^{l+1}\, dv_g     \end{eqnarray*}
It follows from the Sobolev inequality that there exists $C(\omega')>0$ independent of $\eps$ such that
    \begin{eqnarray*} 
    \left ( \int_{\omega'} \left( \eta \ue^{\frac{l+1}{2}} \right)^{\crit } \, dv_g  \right)^{\frac{2}{\crit }} \leq C(\omega') \left(\int_{\omega'} \left| {\nabla} \right( \eta \ue^{\frac{l+1}{2}}    \left)  \right|_g^2 \, dv_g   + B \int_{\omega'} \eta^2 \ue^{l+1} \, dv_g \right)
    \end{eqnarray*}
Combining these inequalities yields
     \begin{equation*}
    \left(  \int_{\omega}      \ue^{\frac{l+1}{2}\crit} \, dv_g  \right)^{\frac{2}{\crit}} \leq   \left(  \int_{\omega'}     \left( \eta \ue^{\frac{l+1}{2}} \right)^{\crit} \, dv_g  \right)^{\frac{2}{\crit}} \leq C\int_{\omega'}\ue^{l+1}\, dv_g
 \end{equation*}
 for $\eps>0$ small enough and where $C$ is independent of $\eps$. Taking $1<l<\crit-1$, we then get that $\ue\to 0$ in $L^q(\omega)$ for some $q>\crit$. Since $\ue$ satisfies \eqref{principal:eq:critical:variete:with:boundary}, it is classical that $\ue\to 0$ in $C^0_{loc}(\omega)$. This proves Claim \ref{claim:2}.\qed
 
  \begin{claim} \label{thm:point:wise:estimate:critical:case} For all $\omega\subset M$ such that $\overline{\omega}\subset M$, there exists $ C(\omega) $  such that
       \begin{equation}\label{estim:weak}
        d_g(x,\xe )^{\frac{n-2}{2}}\ue (x) \leq C(\omega)\hbox{ for all }\eps>0\hbox{ and }x\in\omega.
    \end{equation}
Moreover,
    \begin{equation}\label{estim:weak:bis}
        \lim\limits_{R\to 0} \lim\limits_{{\eps } \to 0} \sup\limits_{x\in \omega \setminus B_{\xe }(R\me )}    d_g(x,\xe )^{\frac{n-2}{2}}|\ue (x) | =0
    \end{equation}
  \end{claim}
\noindent{\it Proof of the Claim:} We argue by contradiction and we let $(y_\eps)_\eps\in \overline{\omega}$ be such that
$$d_g(y_\eps,\xe )^{\frac{n-2}{2}}\ue (y_\eps)=\sup_{x\in\overline{\omega}}d_g(x,\xe )^{\frac{n-2}{2}}\ue (x)\to +\infty\hbox{ as }\eps\to 0.$$
It follows from Claim \ref{claim:2} that $\lim_{\eps\to 0}y_\eps=x_0$. Arguing as in Step 2 of Chapter 4 in Druet-Hebey-Robert \cite{dhr} and using \eqref{estim:out:ball:bis}, we get \eqref{estim:weak}. The second estimate \eqref{estim:weak:bis} follows also from \cite{dhr}.\qed

\medskip\noindent We now state and prove the main result of this section:
\begin{proposition}\label{thm:estimation:forte:u:eps:a:critical:case:article:blowup} Let $\delta>0$ be such that $B_{2\delta}(x_0)\subset M$. Under the assumptions of Theorem \ref{theorem:Mesmar:1:article:blow:up}, we have that
\begin{equation}\label{eq:total}
    \ue(\ye)  =\left(\frac{\me}{\me^2+\frac{f(x_0)}{n(n-2)}d_g(\xe, \ye)^2}\right)^{\frac{n-2}{2}}(1+o(1))+O(\theta_\eps)\hbox{ when }\lim_{\eps\to 0}\ye=x_0.
    \end{equation}
Moreover, there exists $ C(\delta)> 0 $ independent of $ {\mathbb {\epsilon}} $  such that 
    \begin{align}
         \ue (x) & \leq C \frac{\me^{\frac{n-2}{2}  }}{(\me+d_g(x,\xe))^{n-2}} +C \theta_\eps \label{ineq:c0}\\
    |{\nabla} \ue | (x) & \leq C \frac{\me^{\frac{n-2}{2}  }}{(\me+d_g(x,\xe))^{n-1}} +C \theta_\eps \label{ineq:c1}
\end{align}
for all $x\in B_{\delta}(x_0)$ where
\begin{equation}\label{def:theta}
\theta_{\eps } :=\sup_{x \in \partial B_{\delta}(x_0) } \ue (x)\to 0\hbox{ as }\eps\to 0.
\end{equation}
\end{proposition}
\noindent{\it Proof of Proposition \ref{thm:estimation:forte:u:eps:a:critical:case:article:blowup}:} We let $\nu\in (0,1)$ to be fixed later. We let $ \alpha_0> 0 $ such that $ {\Delta_{g}} + \frac {h - \alpha_0} {1 - {\nu}} $ is coercive on $B_{2\delta}(x_0)$ where $ h $ is as in \eqref{lim:h:f}: up to taking $\delta>0$ small, this is always possible.  We let $\tilde {G}_\nu $ be the Green's function of  $ {\Delta_{g}} + \frac {h - \alpha_0} {1 - {\nu_1}} $ on $B_{2\delta} (x_0) $ with Dirichlet boundary condition.  It follows from Robert \cite{robert0existence} that there exists $c_1,c_2>0$ such that 
\begin{equation}\label{control:G}
c_1 d_g(x,y)^{2-n}\leq \tilde{G}_\nu(x,y)\leq c_2 d_g(x,y)^{2-n}\hbox{ for all }x,y\in B_\delta(x_0),\, x\neq y.
\end{equation}
We define the operator 
$$ u\mapsto L_\epsilon u := \Delta_g u+\he u-f \ue^{\crit-2}u,$$
so that \eqref{principal:eq:critical:variete:with:boundary} reads $L_\eps\ue=0$. A straightforward computation yields
   \begin{eqnarray}  \label{401}
        \frac{L_{\eps } {\tilde {G}_\nu^{1-\nu}} }{{\tilde {G}_\nu^{1-\nu}} } ({ x},\xe)& =&  \alpha_0 +\he(x)- h({ x}) +{\nu} ( 1-{\nu} ) \left|\frac{{\nabla} {\tilde {G}_\nu}}{{\tilde {G}_\nu} } \right|_{g}^2 -f\ue^{\crit-2}
  \end{eqnarray}
By standard properties of Green's function \cite {robert0existence}, there exists $  c_1, \rho > 0 $, such that 
\begin{equation} \label{402}
\frac{|{\nabla_{g,x}} {\tilde {G}_\nu}|_{g} }{{\tilde {G}_\nu} } (x,\xe)\geq \frac{c_1}{{d_{g}(x,x_{{\eps }})}} \hbox{ for all }x\in B_{\rho}(\xe)-\{\xe\}.
\end{equation}
Since $ \ue \to 0 $ in $ C_{loc}^0 (M \setminus \{x_0\}) $ and $\he\to h$ in $ C_{loc}^0 (M \setminus \{x_0\}) $,    \eqref{401} yields 
\begin{align*}
     L_{\mathbb {\epsilon}} {\tilde {G}_\nu^{1- \nu}} \geq 0 \hbox{ in }B_{2\delta}(x_0) \setminus B_{\rho}(x_0)
\end{align*}
\medskip\noindent Let $ R> 0 $ to be fixed later. It follows from \eqref{estim:weak:bis} that
\begin{align*}
     {d_{g}(x,x_{{\eps }})}^{2} \ue^{\crit-2} ({ x}) \leq \eta(R)  \hbox{ for all } {x} \in B ({x_{{\mathbb {\epsilon}}}}, \rho) \setminus B ({x_{{\mathbb {\epsilon}}}}, R \me ),
\end{align*}
where $\lim_{R\to +\infty}\eta(R)=0$. Now, using $\he\to h$ in $C^0(M)$, \eqref {401} and \eqref {402}, for any $x\in B(\xe, \rho ) \setminus B(\xe , R \me )$, we get that
\begin{align*}
\frac{L_{\eps } {\tilde {G}_\nu^{1-\nu}}}{ {\tilde {G}_\nu^{1-\nu}}} ({ x},\xe ) & \geq \frac{\alpha_0}{2} + {\nu} (1-{\nu} ) \frac{c_1^2}{{d_{g}(x,x_{{\eps }})}^2}-f \ue^{\crit-2} ({ x}) \\
&\geq \frac{\alpha_0}{2}  + \frac{{\nu} (1-{\nu}) c_1^2 - f \eta(R)}{{d_{g}(x,x_{{\eps }})}^2}\geq 0  
\end{align*}
for $R>0$ large enough. Therefore, we get that
\begin{align}\label{L:eps:G:1:postif:on:ball:r:mu:eps}
    L_{\mathbb {\epsilon}} {\tilde {G}_\nu^{1- \nu}} ({x}, {x_{{\mathbb {\epsilon} }}}) \geq 0\hbox{ for all }x\in B_{\delta}(x_0) \setminus B_{R\me}(\xe).
\end{align}
We fix $\nu_1\in (0,1)$. It follows from \eqref{control:G} and  ${\mathbb {\vert \vert}} {u_{\mathbb {\epsilon}}} {\mathbb {\vert \vert}}_\infty = \me ^{1-n/2}$, that there exists $c_3>0$ such that
\begin{equation} \label{u:eps:leq:G:1:on:ball:r:mu}
    \ue(x)\leq c_3 \me ^{\frac {n-2} {2} - {\nu_1}  (n-2)  }  {\tilde {G}_{\nu_1}^{1-\nu_1}} ({ x},\xe)\hbox{ for all }x \in  \partial B_{ R \me }(\xe).  
\end{equation}
We set $\theta_{\eps } :=\sup_{x \in \partial B_{\delta}(x_0) } \ue (x)$. It follows from Claim \ref{claim:2} that $\lim_{\eps\to 0}\theta_\eps=0$. We fix $\nu_2\in (0,1)$ and we consider the Green's function $\tilde{G}_{\nu_2}$. It follows from \eqref{control:G} that there exists $c_4>0$ such that
\begin{align}\label{u:eps:leq:G:2:on:partial:N:prime}
    \ue ({ x}) \leq  c_4 { {\theta}_{\eps } } {\tilde {G}_{\nu_2}^{1-\nu_2}} (x,\xe)\hbox{ for all }x\in \partial B_{\delta}(x_0).  
\end{align}
We define
$$H_\eps(x):=c_3 \me ^{\frac {n-2} {2} - {\nu_1}  (n-2)  }  {\tilde {G}_{\nu_1}^{1- \nu_1}}(x,\xe) +  c_4 { {\theta}_{\eps } } {\tilde {G}_{\nu_2}^{1-\nu_2}} (x,\xe)\hbox{ for }x\in B_{2\delta}(\xe)-\{\xe\}.$$
It follows from \eqref{L:eps:G:1:postif:on:ball:r:mu:eps}, \eqref{u:eps:leq:G:1:on:ball:r:mu} and \eqref{u:eps:leq:G:2:on:partial:N:prime} that
\begin{equation*}
\left\{\begin{array}{cc}
L_\eps\ue=0\leq L_\eps H_\eps & \hbox{ in }B_\delta(x_0)\setminus B_{ R \me }(\xe))  \\
0<\ue\leq   H_\eps & \hbox{ on }\partial\left(B_\delta(x_0)\setminus B_{ R \me }(\xe)\right)
\end{array}\right.
\end{equation*}
Since $L_\eps\ue\geq 0$ in $\overline{B_\delta(x_0)\setminus B_{ R \me }(\xe)}$, it follows from  \cite {berestycki1994principal0principe0de0max} that
\begin{equation*}
    \ue   \leq H_\eps \hbox{ in }B_\delta(x_0)\setminus B_{R \me } ({x_{{\mathbb {\epsilon}}}}) .
\end{equation*}
Using the pointwise control \eqref{control:G} and that $\Vert\ue\Vert_\infty=\me^{1-n/2}$, we get that for all $\nu_1,\nu_2\in (0,1)$, there exists $C_{\nu_1,\nu_2}>0$ such that
\begin{equation}\label{ineq:nu}
    \ue({ x})  \leq C_{\nu_1,\nu_2}\left( \frac{\me^{\frac{n-2}{2}-{\nu_1}(n-2) }}{(\me+d_g(x,\xe))^{(n-2)(1-\nu_1)}} +\theta_\eps {d_{g}(x,x_{{\eps }})}^{(2-n)(1-{\nu_2}) }  \right )
\end{equation}
for all $x\in B_\delta(x_0)$. Our next step is to prove \eqref{eq:total}. We let $(\ye)_\eps\in M$ such that $\lim_{\eps\to 0}\ye=x_0$.

\smallskip\noindent We first assume that $d_g(\xe,\ye)=O(\me)$ as $\eps\to 0$. Then, \eqref{eq:total} is a direct consequence of \eqref{cv:we:w}. 

\smallskip\noindent From now on, we assume that
\begin{equation*}
\lim_{\eps\to 0} d_g(\xe,\ye)=0\hbox{ and }\lim_{\eps\to 0}\frac{d_g(\xe,\ye)}{\me}=+\infty.
\end{equation*}
We let $G_\eps$ be the Green's function for $\Delta_g+\he$ in $B_{\delta}(x_0)$ with Dirichlet boundary condition. We let $(\ye)_\eps\in B_{\delta/2}(x_0)$. Green's representation formula yields
\begin{equation}\label{green:rep}
    \ue (\ye) = \int_{B_\delta(x_0)} G_\eps(\ye, x) (\Delta_g \ue+\he\ue)(x)\, dv_g(x)-\int_{\partial B_\delta(x_0)} \partial_{\vec{n}} G_\eps(\ye, z)\ue(z)\, d\sigma_g(z).
 \end{equation}
It follows from Robert \cite{robert0existence} that there exists $c_5,c_6>0$ such that 
\begin{equation}\label{est:green}
d_g(x,y)^{n-2}|G_\eps(x,y)|+d_g(x,y)^{n-1}|\nabla_x G_\eps(x,y)|\leq c_5\hbox{ for all }x,y\in B_\delta(x_0),\, x\neq y
\end{equation}
for all $\eps>0$. Combining these estimates with equation \eqref{principal:eq:critical:variete:with:boundary} and \eqref{def:theta}, we get that
\begin{eqnarray}\label{eq:57}
    \ue (\ye)&=& \int_{B_{R\me}(\xe)} G_\eps(\ye, x) f(x)\ue^{\crit-1}(x)\, dv_g(x)+A_\eps(R)+B_\eps
    \end{eqnarray}
where
$$|A_\eps(R)|\leq C \int_{B_\delta(x_0)\setminus B_{R\me}(\xe)} d_g(x,\ye)^{2-n}\ue^{\crit-1}(x)\, dv_g(x)$$
$$\hbox{ and }|B_\eps|\leq    C \int_{\partial B_\delta(x_0)} d_g(z,\ye)^{1-n}\ue(z)\, d\sigma_g(z)\leq C\theta_\eps.$$ 
We deal with the first term of \eqref{eq:57}. With a change of variable and \eqref{cv:we:w}, we get that
\begin{eqnarray*}
&&\int_{B_{R\me}(\xe)} G_\eps(\ye, x) f(x)\ue^{\crit-1}(x)\, dv_g(x)\\
&&= \me^{n-2}\int_{B_{R}(0)} G_\eps(\ye, \hbox{exp}_{\xe}(\me X)) f(\hbox{exp}_{\xe}(\me X))\we^{\crit-1}(X)\, dv_{g_\eps}(X)
\end{eqnarray*}
It follow from \cite{robert0existence} that for any $(z_\eps)_\eps\in M$ such that $\lim_{\eps\to 0}d_g(z_\eps,\xe)=0$ we have that
$$\lim_{\eps\to 0}d_g(\xe, z_\eps)^{n-2}G_\eps(\xe, z_\eps)=\frac{1}{(n-2)\omega_{n-1}}.$$
Since $\me=o(d_g(\xe, \ye))$, we then get that
\begin{eqnarray}
&&\int_{B_{R\me}(\xe)} G_\eps(\ye, x) f(x)\ue^{\crit-1}(x)\, dv_g(x)\nonumber\\
&&= \frac{f(x_0) \me^{n-2}}{(n-2)\omega_{n-1}d_g(\xe,\ye)^{n-2}}\left(\int_{B_{R}(0)} w^{\crit-1}(X)\, dX+o(1)\right)\nonumber\\
&&= \frac{f(x_0) \me^{n-2}}{(n-2)\omega_{n-1}d_g(\xe,\ye)^{n-2}}\left(\int_{\rn}  w^{\crit-1}(X)\, dX+o(1)+\eta(R)\right)\label{eq:78}
\end{eqnarray}
where $\lim_{R\to +\infty}\eta(R)=0$. With \eqref{equ:w:critical:case} and \eqref{cv:we:w}, we get that
\begin{eqnarray}
f(x_0)\int_{\rn}  w^{\crit-1}(X)\, dX&=&\lim_{R\to +\infty}\int_{B(0,R)}\Delta w\, dX=\lim_{R\to +\infty}\int_{\partial B(0,R)}(-\partial_\nu w)\, d\sigma\nonumber\\
&=&\left(\frac{n(n-2)}{f(x_0)}\right)^\frac{n-2}{2}(n-2)\omega_{n-1}.\label{eq:79}
\end{eqnarray}

\medskip\noindent We now deal with $A_\eps(R)$. Using the pointwise control \eqref{ineq:nu}, we get that
\begin{eqnarray}
    |A_\eps(R)|&\leq & C_{\nu_1,\nu_2}\int_{B_\delta(x_0)\setminus B_{R\me}(\xe)} d_g(x,\ye)^{2-n}   \frac{\me^{\frac{n+2}{2}-{\nu_1}(n+2) }}{(\me+d_g(x,\xe))^{(n+2)(1-\nu_1)}} \, dv_g(x)\nonumber\\
    &&+ C_{\nu_1,\nu_2} \theta_\eps^{\crit-1}\int_{B_\delta(x_0)} d_g(x,\ye)^{2-n}  {d_{g}(x,x_{{\eps }})}^{-(n+2)(1-{\nu_2}) }\, dv_g(x) +C    \theta_\eps.\label{remain:34}
 \end{eqnarray}
It follows from Giraud's lemma (see Appendix A of \cite{dhr} for instance) that for $\alpha,\beta\in (0,n-2)$ such that $\alpha+\beta>n$, there exists $C>0$ such that 
$$\int_{B_\delta(x_0)}d_g(y, x)^{\alpha-n}d_g(x,y)^{\beta-n}\, dv_g(x)\leq C\hbox{ for all }y,z\in B_\delta(x_0).$$
Taking $1-\nu_2>0$ close to $0$, we then get that
\begin{equation}\label{eq:80}
\theta_\eps^{\crit-1}\int_{B_\delta(x_0)} d_g(x,\ye)^{2-n}  {d_{g}(x,x_{{\eps }})}^{-(n+2)(1-{\nu_2}) }\, dv_g(x)\leq C\theta_\eps^{\crit-1}\leq C\theta_\eps.
\end{equation}
We now deal with the remaining term of \eqref{remain:34}.  We split the domain $B_\delta(x_0)=D_\eps^1\cup D^2_\eps$ where
$$D_\eps^1:=\{x\in B_\delta(x_0)\hbox{ s.t. }d_g(x,\ye)\geq d_g(\xe,\ye)/2\}$$
$$\hbox{ and }D_\eps^2:=\{x\in B_\delta(x_0)\hbox{ s.t. }d_g(x,\ye)<d_g(\xe,\ye)/2\}.$$
We fix $R>0$. With the change of variable $x:=\hbox{exp}_g(\me X)$, we get that 
\begin{eqnarray}
&&\int_{D_\eps^1\setminus B_{R\me}(\xe)} d_g(x,\ye)^{2-n}   \frac{\me^{\frac{n+2}{2}-{\nu_1}(n+2) }}{(\me+d_g(x,\xe))^{(n+2)(1-\nu_1)}} \, dv_g(x)\nonumber\\
&&\leq C d_g(\xe,\ye)^{2-n}\int_{B_{2\delta}(\xe)\setminus B_{R\me}(\xe)}  \frac{\me^{\frac{n+2}{2}-{\nu_1}(n+2) }}{(\me+d_g(x,\xe))^{(n+2)(1-\nu_1)}} \, dv_g(x)\nonumber\\
&&\leq C d_g(\xe,\ye)^{2-n}\me^{\frac{n-2}{2}}\int_{\rn\setminus B_{R}(0)}  \frac{1}{(1+|X|)^{(n+2)(1-\nu_1)}} \, dX\nonumber\\
&&\leq \eta(R) d_g(\xe,\ye)^{2-n}\me^{\frac{n-2}{2}}\hbox{ where }\lim_{R\to +\infty}\eta(R)=0\label{eq:68}
\end{eqnarray}
when $\nu_1<2/(n+2)$. Concerning the other integral, note that for all $x\in D_\eps^2$, we have that $d_g(x,\xe)\geq d_g(\xe,\ye)/2$. Therefore
\begin{eqnarray}
&&\int_{D_\eps^2} d_g(x,\ye)^{2-n}   \frac{\me^{\frac{n+2}{2}-{\nu_1}(n+2) }}{(\me+d_g(x,\xe))^{(n+2)(1-\nu_1)}} \, dv_g(x)\nonumber\\
&&\leq C \frac{\me^{\frac{n+2}{2}-{\nu_1}(n+2) }}{d_g(\ye,\xe)^{(n+2)(1-\nu_1)}} \int_{d_g(x,\ye)<d_g(\xe,\ye)/2} d_g(x,\ye)^{2-n}\, dv_g(x)\nonumber\\
&&\leq C \frac{\me^{\frac{n+2}{2}-{\nu_1}(n+2) }}{d_g(\ye,\xe)^{(n+2)(1-\nu_1)}}  d_g(\xe,\ye)^2\nonumber\\
&&\leq C \frac{\me^{\frac{n-2}{2}}}{d_g(\xe,\ye)^{n-2}}\left(\frac{\me }{d_g(\ye,\xe)}\right)^{2-\nu_1(n+2)} =o\left(\frac{\me^{\frac{n-2}{2}}}{d_g(\xe,\ye)^{n-2}}\right)\label{eq:69}
\end{eqnarray}
 when $\nu_1<2/(n+2)$. Putting \eqref{eq:78}, \eqref{eq:79}, \eqref{remain:34}, \eqref{eq:80}, \eqref{eq:68} and \eqref{eq:69} together yields
\begin{equation*}
\ue(\ye)=\left(\frac{n(n-2)}{f(x_0)}\right)^\frac{n-2}{2}\frac{\me^{\frac{n-2}{2}}}{d_g(\xe,\ye)^{n-2}}(1+o(1))+O(\theta_\eps)
\end{equation*}
This yields \eqref{eq:total} since $d_g(\xe,\ye)/\me\to +\infty$ as $\eps\to 0$.

\smallskip\noindent When $d_g(\xe,\ye)=o(1)$, \eqref{ineq:c0} is a direct consequence of \eqref{eq:total}. Since $\ue\to 0$ in $C^0_{loc}(M-\{x_0\})$ and $\ue>0$, it follows from Harnack's inequality that there exists $c(\tau)>0$ such that $\ue(x)\leq c(\tau)\ue(y)$ for all $x,y\in B_{2\delta}(x_0)\setminus B_{\tau}(x_0)$. Therefore, if $(\ye)_\eps\in B_{2\delta}(x_0)$ is such that $\ye\not\to x_0$, we have that $\ue(\ye)=O(\theta_\eps)$. This proves \eqref{ineq:c0} when $x$ is far from $x_0$. This proves  \eqref{eq:total} holds in all cases, which yields \eqref{ineq:c0}. 

\smallskip\noindent Concerning the gradient estimate, differentiate Green's representation formula \eqref{green:rep} to obtain
\begin{equation*}
   \nabla \ue (\ye) = \int_{B_\delta(x_0)} \nabla_y G_\eps(\ye, x) (\Delta_g \ue+\he\ue)(x)\, dv_g(x)-\int_{\partial B_\delta(x_0)} \partial_{\vec{n}} \nabla_yG_\eps(\ye, z)\ue(z)\, d\sigma_g(z).
 \end{equation*}
Using the pointwise control \eqref{est:green}, we then get that
\begin{equation*}
|\nabla \ue (\ye)|\leq C \int_{B_\delta(x_0)} d_g(\ye, x)^{1-n}\ue^{\crit-1}(x)\, dv_g(x)+C\int_{\partial B_\delta(x_0)} d_g(\ye, z)^{-n}\ue(z)\, d\sigma_g(z).
 \end{equation*}
We  get \eqref{ineq:c1} arguing as in the proof of \eqref{ineq:c0}. This proves  Proposition \ref{thm:estimation:forte:u:eps:a:critical:case:article:blowup}.\qed

\section{Speed of convergence of $(\xe)_\eps$}
Let $ \Omega \subset {\mathbb {R}}^n $ be a smooth bounded domain. Let $ u \in C^2 (\bar {\Omega}) $, $u>0$,  and $ f \in C^1 (\bar {\Omega}) $ be functions and $ c \in {\mathbb {R}} $. Then for all $z\in \rn$, the Pohozaev identity writes 
\begin{eqnarray}\label{poho:1}
  &&  \int_{\Omega}  \left( (x-z)^i \partial_i u + \frac{n-2}{2}u \right)\left( \Delta_\xi u - c f u^{\crit-1}  \right) dx\\
      &&= \int_{\partial \Omega } \left[ (x-z,\nu) \left( \frac{|\nabla u|_\xi^2}{2} - \frac{c f  u^{\crit}}{\crit} \right) - \left( (x-z)^i \partial_i u + \frac{n-2}{2}u \right) \partial_{\nu} u  \right] \mathrm{d}  \sigma\nonumber\\
      &&  + \frac{1}{\crit} \int_{\Omega} c \langle {\nabla} f(x), x-z\rangle_\xi u^{\crit} dx \nonumber
\end{eqnarray}
Differentiating with respect to $z$, we get that for any $ j \in \{ 1, ..., n \} $ 
  \begin{eqnarray}\label{poho:2}
  &&  - \int_{\Omega}    \partial_j u \left( \Delta_\xi u - c f u^{\crit-1}  \right) \mathrm{dx}  \\        &&= \int_{\partial \Omega } \left[ - \nu_j \left( \frac{|\nabla u|_\xi^2}{2} - \frac{c f  u^{\crit}}{\crit}  \right)  +  \partial_j u \,   \partial_{\nu} u \right]  \mathrm{  d}\sigma 
       - \frac{c}{\crit} \int_{\Omega} { \partial_j} f(x)  u^{\crit} \mathrm{ d}x \nonumber
\end{eqnarray} 
We refer to  Ghoussoub-Robert  \cite{gmr} for a proof. We fix $\delta\in (0, i_g(M,x_0))$. We define
\begin{equation*}
\hue(X):=\ue \left( \exp_{\xe} (X) \right) \hbox{ for all } X \in B_{\delta}(0) \subset {\mathbb{R}}^n .
\end{equation*}
Therefore, equation \eqref{principal:eq:critical:variete:with:boundary} rewrites
\begin{equation*}
\Delta_{\hge}\hue+\hue\ue=\hfe\hue^{\crit-1}\hbox{ in }B_{\delta}(0).
\end{equation*} 
where $\hhe(X):=\he \left( \exp_{\xe} (X) \right)$ and $\hfe(X):=f \left( \exp_{\xe} (X) \right)$ for all $X \in B_{\delta}(0) \subset {\mathbb{R}}^n$
and $\hge:=\hbox{exp}^\star_{\xe }g$ is the pull-back of $g$ via the exponential map.

\begin{lemma} \label{lem:cv:lq}
 Let $ (\phi_\eps)_\eps \in C^0 (B_\delta (0)) $  such that 
 $$\left\{\begin{array}{c}
 \lim_{\eps\to 0}\phi_{\mathbb {\epsilon}} (0) =s\in\rr,\\
 |\phi_\eps(X)-\phi_\eps(0)|\leq C|X|\hbox{ for all }X\in B_\delta(0)\hbox{ and }\eps>0.
 \end{array}\right.$$
 We fix $p\geq 0$ and $q\geq 1$. Then
 \begin{eqnarray*}
 && \int_{B_\delta(0)} \phi_{\eps } |X|^p\hue^q \, dX \\
 &&=\left\{\begin{array}{cc}
 \me^{n+p - \frac{q(n-2)}{2}}  \left(  s \int_{{\mathbb{R}}^n } |X|^p{w }^q \mathrm{d}X + o(1) \right)+O(\theta_\eps^q)&\hbox{ if }(n-2)q>p+n,\\
 \left(s\left(\frac{n(n-2)}{f(x_0)}\right)^{q\frac{n-2}{2}}\omega_{n-1}+o(1)\right)\me^{\frac{q(n-2)}{2}} \ln\left(\frac{1}{\me}\right)+O(\theta_\eps^q)&\hbox{ if }(n-2)q=p+n,\\
 O(\me^{q\frac{n-2}{2}})+O(\theta_\eps^q)&\hbox{ if }(n-2)q<p+n
 \end{array}\right.\end{eqnarray*}
Moreover, for any family $(\delta_\eps)_{\eps}\in (0,1)$ such that $\lim_{\eps\to 0}\delta_\eps=\lim_{\eps\to 0}\frac{\me}{\delta_\eps}=0$, we have that
\begin{equation}\label{lim:addition}
\int_{B_{\delta_\eps}(0)} \phi_{\eps } |X|^p\hue^q \, dX=\me^{n+p - \frac{q(n-2)}{2}}  \left(  s \int_{{\mathbb{R}}^n } |X|^p{w }^q \mathrm{d}X + o(1) \right)+O(\theta_\eps^q)\hbox{ if }q>\frac{p+n}{n-2}.
\end{equation}
\end{lemma}
\noindent{\it Proof:} We fix $\nu>0$. It follows from \eqref{eq:total} and \eqref{ineq:c0} that there exists $\alpha\in (0,\delta)$ such that
$$\left|\hue^q(X)-\left(\frac{\me}{\me^2+\frac{f(x_0)}{n(n-2)}|X|^2}\right)^{q\frac{n-2}{2}}\right|\leq \nu \left(\frac{\me}{\me^2+\frac{f(x_0)}{n(n-2)}|X|^2}\right)^{q\frac{n-2}{2}}+C\theta_\eps^q$$
for all $X\in B_\rho(0)\setminus B_\alpha(0)$. Note that for all $\alpha\in (0,\delta)$, it follows from the Harnack inequality that
\begin{equation*}
\int_{B_\delta(0)\setminus B_\alpha(0)} \phi_{\eps } |X|^p\hue^q \, dX =O(\theta_\eps^q).
\end{equation*}
We then get that
\begin{eqnarray*}
&&\left|\int_{B_\delta(0) } \phi_{\eps } |X|^p\hue^q \, dX -\int_{B_\alpha(0) } \phi_{\eps } |X|^p\left(\frac{\me}{\me^2+\frac{f(x_0)}{n(n-2)}|X|^2}\right)^{q\frac{n-2}{2}} \, dX\right|\\
&&\leq C\nu \int_{B_\delta(0) } |X|^p\left(\frac{\me}{\me^2+|X|^2}\right)^{q\frac{n-2}{2}} \, dX+C\theta_\eps^q
\end{eqnarray*}
We the get Lemma \ref{lem:cv:lq} when $q(n-2)<n+p$. With the change of variable $X=\me Y$, we get that 
\begin{eqnarray*}
&&\int_{B_\alpha(0) } \phi_{\eps } |X|^p\left(\frac{\me}{\me^2+\frac{f(x_0)}{n(n-2)}|X|^2}\right)^{q\frac{n-2}{2}} \, dX\\
&&=\me^{n+p-q\frac{n-2}{2}}\int_{B_{\alpha/\me}(0) } \phi_{\eps }(\me Y) |Y|^p\left(\frac{1}{1+\frac{f(x_0)}{n(n-2)}|Y|^2}\right)^{q\frac{n-2}{2}} \, dY\\
&&=\me^{n+p-q\frac{n-2}{2}}\left\{\begin{array}{cc}
s\int_{\rn}|Y|^pw^q(Y)\, dY+o(1)&\hbox{ if }q(n-2)>p+n\\
s\left(\frac{n(n-2)}{f(x_0)}\right)^{q\frac{n-2}{2}}\omega_{n-1}\ln\left(\frac{1}{\me}\right)+o(\ln\me)&\hbox{ if }q(n-2)=p+n\end{array}\right.
\end{eqnarray*}
The proof of \eqref{lim:addition} is similar by taking $\alpha:=\delta_\eps$. Putting these estimates together yields Lemma \ref{lem:cv:lq}.\qed

\smallskip\noindent We now prove \eqref{d:est:petit:o:de:mu:eps}. We fix $l\in \{1,...,n\}$. We define
$$\delta_\eps:=\me^{\frac{1}{n-1}}.$$
Pohozaev's identity  \eqref{poho:2} applied to  $ \hue $ reads
\begin{align} \label{poho:abcd} 
A_{\eps } = - B_{\eps } +  C_{\eps } - D_{\eps }
\end{align}
With
\begin{align*}
 A_{\eps }& :  = - \frac{1}{\crit} \int_{{B_{\delta_\eps}(0)}} { \partial_l} {{\hat{f}}_{\eps }}  {{\hat{u}}_{{\eps }}}^{\crit} \mathrm{d}X \\
    B_{\eps } & :=  \int_{\partial {B_{\delta_\eps}(0)}}- \left[  \nu_l \left(  \frac{|{\nabla} {{\hat{u}}_{{\eps }}}|^2 }{2} - \frac{{\hat{f}}_\epsilon {\hat{u}}_{\epsilon}^{\crit}}{\crit} \right) +    \partial_l {{\hat{u}}_{{\eps }}}  \partial_\nu {{\hat{u}}_{{\eps }}} \right] \mathrm{d}\nu   \\
    C_{\eps } & :=   \int_{{B_{\delta_\eps}(0)}}      \partial_l {{\hat{u}}_{{\eps }}}  {{\hat{h}}_{{\eps }}} {{\hat{u}}_{{\eps }}} \mathrm{d}X \hbox{ and }D_{\eps }  :=  \int_{{B_{\delta_\eps}(0)}}  \partial_l {{\hat{u}}_{{\eps }}}   (\Delta_{\xi }   {{\hat{u}}_{{\eps }}} - \Delta_{{ {\hat{g}}_{{\eps }}}} {{\hat{u}}_{{\eps }}}  ) \mathrm{d}X 
\end{align*}
We estimate these terms separately. It follows from \eqref{ineq:c0} and \eqref{ineq:c1} that
\begin{equation*}
B_\eps=O\left(\me\left(\left(\frac{\me}{\delta_\eps}\right)^{n-3}+\frac{\me^{n-1}}{\delta_\eps^{n+1}}+\frac{\delta_\eps^{n-1}}{\me}\theta_\eps^2\right)\right)=o(\me)\hbox{ as }\eps\to 0.
 \end{equation*}
Concerning $C_\eps$, integrating by parts, we have that 
 \begin{equation*}
   C_\eps=  \int_{{B_{\delta_\eps}(0)} }  { \partial_l} {{\hat{u}}_{{\eps }}}  {{\hat{h}}_{{\eps }}}  {{\hat{u}}_{{\eps }}} \mathrm{d} X  = -\int_{{B_{\delta_\eps}(0)} }   \partial_l\frac{{{\hat{h}}_{{\eps }}}}{2}   {{\hat{u}}_{{\eps }}}^2 \mathrm{d} X + \int_{\partial {B_{\delta_\eps}(0)}  }  {{\hat{h}}_{{\eps }}} \frac{{{\hat{u}}_{{\eps }}}^2}{2} \vec{\nu}_l\mathrm{d}  \sigma  
\end{equation*}
With \eqref{ineq:c0} and Lemma \ref{lem:cv:lq}, we then get that
\begin{eqnarray*} 
  C_{\eps }   =   O\left(\me\left( o(1)+\left(\frac{\me}{\delta_\eps}\right)^{n-3}+\frac{\delta_\eps^{n-1}}{\me}\theta_\eps^2\right)\right)=o(\me)\hbox{ as }\eps\to 0.
  \end{eqnarray*}
\noindent We now estimate $ D_{\mathbb {\epsilon}} $. We write 
$$ - (\Delta_{{{\hat {g}}_{{\mathbb {\epsilon}}}}} - \Delta_{\xi}) = ({{\hat {g}}_{\mathbb {\epsilon}}^{ij}} - {\delta^{ij}}) {\partial_{ij}} - {{\hat {g}}_{\mathbb {\epsilon}}^{ij}} {\hat {\Gamma}_{ij}^k}(\hge) {\partial_k} $$
where the $ {\hat {\Gamma}_{ij}^k}$'s are the Christoffel symbols of the metric $\hge$. The following lemma is reminiscent in such problems:
\begin{lemma}\label{lem:ipp} Let $\Omega$ be a smooth domain of $\rn$. For any $i,j,k\in \{1,...,n\}$, let us consider $a^{ijk}\in C^1(\rn)$. We assume that $a^{ijk}=a^{jik}$ for all $i,j,k\in\{1,...,n\}$. Then for all $u\in C^2(\rn)$, we have that
\begin{eqnarray*}
\int_\Omega a^{ijk}\partial_{ij}u\partial_k u\, dx&=&\int_\Omega\left(-\partial_la^{lij}+\frac{1}{2}\partial_l a^{ijl}\right)\partial_i u\partial_j u\, dx\\
&&+\int_{\partial\Omega}\left(-\frac{1}{2}a^{ijl}\vec{\nu}_l+ a^{l ji}\vec{\nu}_l\right)\partial_i u\partial_j u \, d\sigma 
\end{eqnarray*}
where $\vec{\nu}$ is the outer normal vector at $\partial\Omega$ and  Einstein's summation convention has been used. 
\end{lemma} 
The proof is by integrations by parts and goes back to Hebey-Vaugon \cite{hv:mz} and can also be found in Cheikh-Ali \cite{hca:pjm}. It follows from Lemma \ref{lem:ipp} that 
\begin{eqnarray*}
\int_{B_{\delta_\eps}(0)} (\hge^{ij}-\delta^{ij})\partial_{ij}\hue\partial_l \hue\, dx&=&\int_{B_{\delta_\eps}(0)}\left(-\partial_m \hge^{mj}\delta_{j,l}+\frac{1}{2}\partial_l \hge^{ij}\delta_{m,l}\right)\partial_i \hue\partial_j\hue\, dx\\
&&+\int_{\partial B_{\delta_\eps}(0)}\left(-\frac{1}{2}(\hge^{ij}-\delta^{ij})\delta_{m,l}\vec{\nu}_m+ (\hge^{mj}-\delta^{mj})\delta_{i,l}\vec{\nu}_l\right)\partial_i \hue\partial_j \hue \, d\sigma 
\end{eqnarray*}
Using  \eqref{ineq:c0} to control the boundary terms, we get that
 \begin{eqnarray*}
 D_{\eps } &=& \int_{B_{\delta_\eps}(0)}\left(-\partial_m \hge^{mj}\delta_{j,l}+\frac{1}{2}\partial_m \hge^{ij}\delta_{m,l}\right)\partial_i \hue\partial_j\hue\, dx\\
 &&-\int_{B_{\delta_\eps}(0)} \hge^{ij}\Gamma_{ij}^k(\hge)\partial_k\hue\partial_l\hue\, dX+ \underbrace{O\left(\me\left( \left(\frac{\me}{\delta_\eps}\right)^{n-3}+\frac{\delta_\eps^{n-1}}{\me}\theta_\eps^2\right)\right)}_{=o(\me)}\end{eqnarray*}
as $\eps\to 0$. Since 
$$\Gamma_{ij}^k(\hge)=\frac{1}{2}\hge^{km}\left(\partial_i(\hge)_{jm}+\partial_j(\hge)_{im}-\partial_m(\hge)_{ij}\right)$$
and $\hge$ is normal at $0$ (that is $\partial_m(\hge)_{ij}(0)=0$ for all $i,j,m\in \{1,...,n\}$), we then get that there exists $a_{ij\alpha}\in\rr$, $i,j,\alpha\in \{1,...,n\}$ such that 
\begin{eqnarray*}
D_\eps&=& \int_{B_{\delta_\eps}(0)}a_{ij\alpha}X^\alpha\partial_i\hue\partial_j\hue\, dx+O\left(\int_{B_{\delta_\eps}(0)}|X|^2|\nabla\hue|^2\, dX\right) + o(\me)\end{eqnarray*}
With \eqref{ineq:c1}, we get that
\begin{eqnarray*}
\int_{B_{\delta_\eps}(0)}|X|^2|\nabla\hue|^2\, dX&\leq& C\int_{B_{\delta_\eps}(0)}|X|^2\frac{\me^{n-2}}{(\me+|X|)^{2(n-1)}}\, dX+C\theta_\eps^2\delta_\eps^n\\
&\leq & C \me^{2}\int_{B_{\delta_\eps/\me}(0)}\frac{|X|^2\, dX}{(1+|X|)^{2(n-1)}}+C\theta_\eps^2\delta_\eps^n=o(\me)
\end{eqnarray*}
With \eqref{ineq:c1}, given $R>0$, using \eqref{control:theta} and $n>3$, we have that
\begin{eqnarray*}
&&\left|\int_{B_{\delta_\eps}(0)\setminus B_{R\me}(0)} X^\alpha\partial_i\hue\partial_j\hue\, dx\right|\leq   C\int_{B_{\delta_\eps}(0)\setminus B_{R\me}(0)} \frac{\me^{n-2}|X|\, dX}{(\me+|X|)^{2(n-1)}}+C\theta_\eps^2\delta_\eps^n\\
&&\leq C\me\int_{\rn\setminus B_R(0)} \frac{|Y|\, dX}{(1+|Y|)^{2(n-1)}}+C\theta_\eps^2\delta_\eps^n\leq \eta(R)\me+o(\me)
\end{eqnarray*}
where $\lim_{R\to +\infty}\eta(R)=0$. Using the change of variable $X=\me Y$, the convergence \eqref{cv:we:w} and the radial symmetry of $w$, we have that
\begin{eqnarray*}
\int_{B_{R\me}(0)} X^\alpha\partial_i\hue\partial_j\hue\, dx&=& \me\int_{B_{R}(0)} Y^\alpha\partial_i\we\partial_j\we\, dY\\
&=& \me \left(\int_{B_R(0)}X^\alpha\partial_i w\partial_jw\, dY+o(1)\right)=o(\me)\hbox{ since }n>3.
\end{eqnarray*}
Therefore, we get that $\int_{B_{\delta_\eps}(0)}X^\alpha\partial_i\hue\partial_j\hue\, dx=o(\me)$, and then
$$D_\eps=o(\me)\hbox{ for }n\geq 4.$$

\medskip\noindent We now deal with $A_\eps$. With a Taylor expansion of $ f$, we get 
\begin{eqnarray*}  A_{\eps } &= & - \frac{1}{\crit} \int_{{B_{\delta_\eps}(0)}} { \partial_l} {{\hat{f}}_{\eps }}  {{\hat{u}}_{{\eps }}}^{\crit} \mathrm{d}X \\
&   =&   - \frac{1}{\crit} { \partial_l} {{\hat{f}}_{\eps }}(0)  \int_{{B_{\delta_\eps}(0)}} {{\hat{u}}_{{\eps }}}^{\crit} \mathrm{d}X - \frac{1}{\crit} { \partial_{lj}} {{\hat{f}}_{\eps }} (0)  \int_{{B_{\delta_\eps}(0)} }  X^j {{\hat{u}}_{{\eps }}}^{\crit} \mathrm{d}X   \\
&   & +O \left(  \int_{{B_{\delta_\eps}(0)}} |X|^2 {{\hat{u}}_{{\eps }}}^{\crit } \mathrm{d}X \right)
\end{eqnarray*}
Arguing as above, we get that 
$$\int_{{B_{\delta_\eps}(0)}} |X|^2 {{\hat{u}}_{{\eps }}}^{\crit } \mathrm{d}X =o(\me)\hbox{ and } \int_{{B_{\delta_\eps}(0)} }  X^j {{\hat{u}}_{{\eps }}}^{\crit} \mathrm{d}X =o(\me)\hbox{ as }\eps\to 0.$$
With \eqref{eq:total}, we get that there exists $C_0>0$ such that
$$\int_{{B_{\delta_\eps}(0)}} {{\hat{u}}_{{\eps }}}^{\crit} \mathrm{d}X=C_0+o(1)\hbox{ as }\eps\to 0.$$
Therefore, we get that
\begin{equation*}
A_\eps= \left(-\frac{C_0}{\crit}+o(1)\right){ \partial_l} {{\hat{f}}_{\eps }}(0)+o(\me)\hbox{ as }\eps\to 0.
\end{equation*}
Putting the estimates of $A_\eps$, $B_\eps$, $C_\eps$ and $D_\eps$ into \eqref{poho:abcd} yields
\begin{equation}
\partial_l\hfe(0)=o(\me)\hbox{ as }\eps\to 0\hbox{ for all }l\in \{1,...,n\}.\label{eq:98}
\end{equation}
Passing to the limit, we get that $\nabla f(x_0)=0$. We now express $\partial_l\hfe(0)$ more precisely. We write
$$ \hfe(X)=f \circ \exp_{\xe}  (X) = \tilde {f} \circ \varphi (X_{\mathbb {\epsilon}}, X) \hbox {for} X \in \mathbb {R}^n $$
where  $ \tilde {f}: = f \circ \exp_{x_0} $ and $ \varphi (Z, X): = \exp_{x_0}^{- 1} \circ \exp_{\exp_{x_0} (Z)} (X) $ for $X,Z\in\rn$. We set $ X_{\mathbb {\epsilon}}: = \exp_{{x_{{0}}}}^{- 1} ({\xe })$. Since $\nabla\tilde{f}(0)=0$, we get that
 \begin{equation*}
 { \partial_l} ( f \circ \exp_{\xe } ) (0) = \frac{\partial^2 (f \circ \exp_{x_0}  )}{ \partial x_l  \partial x_j  } (0)  X_\eps^j + o ( |X_\eps|)
 \end{equation*}
therefore, with \eqref{eq:98}, we get that 
$$ \frac{\partial^2 (f \circ \exp_{x_0}  )}{ \partial x_l  \partial x_j  } (0)  X_\eps^j =o ( |X_\eps|)+o(\me)\hbox{ as }\eps\to 0\hbox{ for all }l=1,..,n.$$
Since $\nabla^2f(x_0)$ is nondegenerate, we then get that $|X_\eps|=o(\me)$, in other words $d_g(\xe, x_0)=o(\me)$ as $\eps\to 0$. This proves \eqref{d:est:petit:o:de:mu:eps}.

\begin{lemma} \label{control:theta:mu}
 Under the assumptions of Theorem \ref{theorem:Mesmar:1:article:blow:up}, we have that
 \begin{eqnarray} \label{control:theta}
    \theta_\eps   = \left\{ \begin{array}{ll}
           o (1 )    & \hbox{ if }  n\geq 7,\\  
           o (\me )    & \hbox{ if }  n\in \{5,6\}.  \\
 o\left( \me \sqrt{\ln( \frac{1}{\me}}\right)   & \hbox{ if } n= 4.
      \end{array} \right.
  \end{eqnarray} 
 \end{lemma}
\noindent{\it Proof:} The case $n\geq 7$ is simply \eqref{def:theta}. It follows from \eqref{ineq:c0} that
\begin{eqnarray*}
\int_{B_\delta(x_0)}\ue^2\, dv_g &\leq & \int_{B_{2\delta}(\xe)}\frac{\me^{n-2}}{(\me+d_g(x,\xe))^{2(n-2)}}\, dv_g +C\theta_\eps^2\\
&\leq & \me^2\int_{B_{2\delta/\me}(0)}\frac{1}{(1+|X|))^{2(n-2)}}\, dv_g +C\theta_\eps^2\\
&\leq & C\left\{ \begin{array}{ll}
           \me^2    & \hbox{ if }  n \geq 5  \\
 \me^2 \ln( \frac{1}{\me})   & \hbox{ if } n= 4  \\
      \end{array} \right.
\end{eqnarray*}
Equation \eqref{principal:eq:critical:variete:with:boundary} rewrites $\Delta_g \ue+(\he-f\ue^{\crit-2})\ue=0$ in $M$. Since $\ue\to 0$ in $C^0_{loc}(M-\{x_0\})$ and $\ue>0$, it follows from Harnack's inequality that there exists $c>0$ such that $\ue(x)\leq c\ue(y)$ for all $x,y\in B_{2\delta}(x_0)\setminus B_{\delta/3}(x_0)$. Therefore, with the definition \eqref{def:theta} of $\theta_\eps$, we get that 
$$\int_{B_{2\delta}(x_0)\setminus B_{\delta}(x_0)}\ue^2\, dv_g\geq c^{-2}\theta_\eps^2.$$
When $4\leq n\leq 6$, it follows from $L^2-$concentration assumption \eqref{hyp:conctration:L2:article:blowup} that
$$\int_{B_{2\delta}(x_0)\setminus B_{\delta}(x_0)}\ue^2\, dv_g\leq \int_{M\setminus B_\delta(x_0)}\ue^2\, dv_g=o\left(\int_{B_\delta(x_0)}\ue^2\, dv_g\right)\hbox{ as }\eps\to 0.$$
Putting these inequalities together, we get \eqref{control:theta}. This proves   Lemma \ref{control:theta:mu}.\qed

\section{Interaction with the scalar curvature: proof of \eqref{relation:de:h0:avec:courbure:scalaire:cas:variete:a:bord}}
This part is strongly inspired by Cheikh-Ali \cite{hca:pjm}. We define
$$\delta_\eps:=\left\{\begin{array}{cc}
\me^{\frac{2}{n-2}} &\hbox{ if }n\geq 7\\
\delta & \hbox{ if }n\in \{4,5,6\}. 
\end{array}\right.$$
Writing the Pohozaev identity \eqref{poho:1} for $\hue$ that satisfies \eqref{principal:eq:critical:variete:with:boundary}, we get that
\begin{align} \label{poho:abcd:2}
    A_{\eps } +  B_{\eps } = C_{\eps } + D_{\eps } 
\end{align}
where
\begin{align*}
    B_{\eps } & := \int_{\partial {B_{\delta_\eps}(0)}} \left[(X,\nu) \left(  \frac{|{\nabla} {{\hat{u}}_{{\eps }}}|_\xi^2 }{2} - \frac{{{\hat{f}}_{\eps }} {{\hat{u}}_{{\eps }}}^{\crit}}{\crit}     \right) - \left(   X^l \partial_l {{\hat{u}}_{{\eps }}} + \frac{n-2}{2} {{\hat{u}}_{{\eps }}}  \right) \partial_\nu {{\hat{u}}_{{\eps }}}\right] \mathrm{d}  \nu   \\
    C_{\eps } & := - \int_{{B_{\delta_\eps}(0)}} \left(    X^l \partial_l {{\hat{u}}_{{\eps }}} + \frac{n-2}{2} {{\hat{u}}_{{\eps }}}   \right) {{\hat{h}}_{{\eps }}} {{\hat{u}}_{{\eps }}} \mathrm{d} X \\
    D_{\eps } & := - \int_{{B_{\delta_\eps}(0)}} \left(      X^l \partial_l {{\hat{u}}_{{\eps }}} + \frac{n-2}{2} {{\hat{u}}_{{\eps }}} \right)  ( \Delta_{{{\hat{g}}_{{\eps }}}} {{\hat{u}}_{{\eps }}} - \Delta_{\xi } {{\hat{u}}_{{\eps }}}  ) \mathrm{d} X \\
   A_{\eps }& :  =  \frac{1}{\crit} \int_{{B_{\delta_\eps}(0)}} ( {\nabla} {{\hat{f}}_{\eps }}, X) {{\hat{u}}_{{\eps }}}^{\crit} \mathrm{d} X
\end{align*}
Following Cheikh-Ali \cite{hca:pjm} and using the pointwise controls \eqref{eq:total}, \eqref{ineq:c0}, \eqref{ineq:c1} and the control \eqref{control:theta} on $(\theta_\eps)_\eps$ when $n\in \{4,5,6\}$, we get that

\begin{equation*}
    B_{\eps } = \left\{\begin{array}{cc}
o(\me^2) &\hbox{ if }n\geq 5\\
o\left(\me^2\ln\frac{1}{\me}\right)  & \hbox{ if }n=4. 
\end{array}\right.\hbox{ as }\eps\to 0,
\end{equation*}
 \begin{eqnarray*}
    C_{\eps }    =   \left\{ \begin{array}{ll}
           h (x_0) \me^2 \ln \left( \frac{1}{\me} \right) \left(\frac{8}{f(x_0)}\right)^{2} \omega_3 + o(\me^2\ln\frac{1}{\me })  & \hbox{ if }  n =4   \\
          h (x_0) \me^2 \int_{{\mathbb{R}}^n} {w }^2 \mathrm{d}X + o(\me^2) &  \hbox{ if } n \geq 5   
      \end{array} \right.
  \end{eqnarray*}
 \begin{eqnarray*} 
    D_{\eps }     =   \left\{ \begin{array}{ll}
           - \me^2 \ln \left( \frac{1}{\me} \right)  \frac{1}{6} \hbox{Scal}_{g} (\xe) \left(\frac{8}{f(x_0)}\right)^{2} \omega_3 + o(\me^2\ln\frac{1}{\me })  & \hbox{ if }  n =4   \\
         - \me^2  \frac{n-2}{4(n-1)} \hbox{Scal}_{g} (\xe) \int_{{\mathbb{R}}^n} {w }^2 \mathrm{d} X    +   o(\me^2)  &  \hbox{ if } n \geq 5   
      \end{array} \right.
      \end{eqnarray*}
We are then left with estimating $A_\eps$. With a Taylor expansion of $\hfe$, we get that
 \begin{eqnarray*}
   A_{\eps } &= & \frac{1}{\crit}     \int_{{B_{\delta_\eps}(0)} }  { \partial_i} {{\hat{f}}_{\eps }}(0)  X^i   {{\hat{u}}_{{\eps }}}^{\crit} \mathrm{d} X     + \frac{1}{\crit} { \partial_{ij}} {{\hat{f}}_{\eps }} (0)  \int_{{B_{\delta_\eps}(0)} } X^i X^j {{\hat{u}}_{{\eps }}}^{\crit} \mathrm{d} X \\
 &&  + O \left(  \int_{{B_{\delta_\eps}(0)}} |X|^3 {{\hat{u}}_{{\eps }}}^{\crit } \mathrm{d} X  \right)
\end{eqnarray*}
Since $\hfe:=f \circ \hbox{exp}_{\xe}$ and $\nabla f(x_0)=0$, we get that $\nabla\hfe(0)=O(d_g(\xe, x_0))$. With \eqref{d:est:petit:o:de:mu:eps}, we then get that $\nabla\hfe(0)=o(\me)$. It follows from  Lemma \ref{lem:cv:lq} that $ \int_{{B_{\delta_\eps} (0)}} | X |^3 {{\hat {u}}_{{\mathbb {\epsilon}}}}^{\crit} \mathrm {d} X = o (\me^2) $ and $ \int_{ {B_{\delta_\eps} (0)}} | X | {{\hat {u}}_{{\mathbb {\epsilon}}}}^{\crit} \mathrm {d} X = O ({\mu_{\mathbb {\epsilon} }}) $. Therefore, we get that

\begin{equation*}
 A_{\eps }=  \frac{1}{\crit} { \partial_{ij}} {{\hat{f}}_{\eps }} (0)  \int_{{B_{\delta_\eps}(0)} } X^i X^j {{\hat{u}}_{{\eps }}}^{\crit} \mathrm{d} X+ o(\me^2) 
 \end{equation*}
Arguing as in the proof of Lemma \ref{lem:cv:lq}, we get
\begin{eqnarray*}   \int_{{B_{\frac{\delta_\eps}{\me  }}(0)}} X^i X^j {{\hat{u}}_{{\eps }}}^{\crit} \mathrm{d} X    =  \me^2  \int_{{\mathbb{R}}^n} X^i X^j   w^{\crit} \mathrm{d} X  + o( \me ^2)  \hbox{ when } n \geq 4.
\end{eqnarray*} 
Since $w$ is radially symmetric, we get that  $\int_{{\mathbb{R}}^n} X^i X^j   w^{\crit} \mathrm{d} X=\frac{\delta_{ij}}{n}\int_{{\mathbb{R}}^n} |X|^2   w^{\crit} \mathrm{d} X$. Since $\hge$ is normal at $0$, we have that $\Delta_g f(\xe)=-\sum_i\partial_{ii}\hfe(0)$, which yields
We claim that
 \begin{eqnarray*}
 A_{\eps }   =   - \frac{1}{ \crit n} \Delta_{g} f(x_0)  \me^2  \int_{{\mathbb{R}}^n} |X|^2 w^{\crit} \mathrm{d} X  +o(\me^2) & \hbox{ if } n \geq 4
 \end{eqnarray*} 
By Jaber \cite {jaber2015optimal} we have that  
$$\frac{\int_{{\mathbb{R}}^n} |X|^2 w^{\crit} \mathrm{d} X }{ \int_{{\mathbb{R}}^n}w^2 \mathrm{d} X } = \frac{n^2(n-4)}{4(n-1)f(x_0)}\hbox{ for }n\geq 5.$$
Putting the expressions of $A_\eps$, $B_\eps$, $C_\eps$ and $D_\eps$ in \eqref{poho:abcd:2} and letting $\eps\to 0$ yield \eqref{relation:de:h0:avec:courbure:scalaire:cas:variete:a:bord}. This ends the proof of Theorem \ref{theorem:Mesmar:1:article:blow:up}.

\section{Application to a super-critical problem: proof of Theorem \ref{th:mesmar:2:supercitical:case:blow:up} }
We follow the notations and assumptions of Theorem \ref{th:mesmar:2:supercitical:case:blow:up}. We consider a family $ (u_ \epsilon)_{\epsilon> 0} \in C^2_G (X) $ of $ G-$invariant solutions to the problem
\begin{equation} \label{main:eq}
    \Delta_g\ue + \he \ue ={{\mathbb{\lambda}}_{\eps }}\ue^{\critk - 1} \, ,\,  \int_X\ue^{\critk}{\mathrm dv_g} = 1\, ,\, \Vert\ue\Vert_2\to 0\hbox{ as }\eps\to 0
\end{equation}
 where $ (h_ \epsilon)_{\epsilon> 0} \in C^1_G(X) $  is such that there exists $ h  \in C^1_G(X)$  such that \eqref{heps:intro:eng} holds and $(\le)_\eps$ is such that \eqref{align:hyp:hme:blowup:supercritical} holds. 

\begin{claim} There exists $x_0\in X$ such that 
\begin{equation}\label{eq:conc}
\lim_{\eps\to 0}\int_{B_{\delta}(Gx_0)}\ue^{\critk}\, dv_g=1\hbox{ for all }\delta>0.
\end{equation}
\end{claim}
\noindent{\it Proof:} We fix a point $z_0\in X$. We choose $ \eta_0 \in C^{\infty} ({\mathbb {R}}) $ such that $\eta_0(t)=1$ for $t\leq 1$ and $\eta_0(t)=0$ for $t\geq 2$.  Given $\delta>0$, we define $\eta (x) := \eta_0 (\frac {d_{ {g }} (Gx_0, x)} {\delta}) $ for all $x\in X$. For $\delta>0$ small enough, we have that $\eta\in C^\infty_G(X)$.  Multiply \eqref{principale:equation:supercritical:blowup} by $ \eta^2 \ue ^l $ for some $ l> 1 $ and integrate over $X$, we get that
 \begin {equation} \label{equality:with:eta2:and:u:eps:l}
     \int_M \eta^2 \ue ^l \Delta_g \ue  {\mathrm dv_g} + \int_M \eta^2 {h_{\mathbb { \epsilon}}} \ue ^{l + 1} {\mathrm dv_g} = {{\mathbb {\lambda}}_{\mathbb {\epsilon}}} \int_M \eta^2 \ue ^{l + \critk -1} {\mathrm dv_g}
 \end {equation}
 As in the proof of Claim 2, we get that
\begin{eqnarray*}
&& \int_M \left| {\nabla}  \right( \eta \ue^{\frac{l+1}{2}}    \left)  \right|_g^2 \, dv_g  = \frac{(l+1)^2}{4l} \int_M \eta^2 \ue^l \Delta_g \ue \, dv_g  + \frac{l+1}{2l}\int_M \left( |{\nabla}  \eta |_g^2 + \frac{l-1}{l+1} \eta \Delta_g \eta \right) \ue^{l+1} \, dv_g  
 \end{eqnarray*}
Using \eqref{equality:with:eta2:and:u:eps:l} and H\"older's inequality, we get
    \begin{eqnarray} \label{ineq:eta:2:u:alpha}
  &&  \int_X \left| {\nabla}  \right( \eta \ue^{\frac{l+1}{2}}    \left)  \right|_g^2 \, dv_g\leq C \int_X \ue^{l+1} \, dv_g    \\
  && +    \frac{(l+1)^2}{4l}\le  \left ( \int_X \left( \eta \ue^{\frac{l+1}{2}} \right)^{\critk} \, dv_g  \right)^{\frac{2}{\critk}} \left( \int_{B_{2\delta}(Gz_0)}  \ue^{\critk} \, dv_g   \right)^{1-\frac{2}{\critk}} \nonumber
    \end{eqnarray}
It follows from Faget \cite{faget2002best} that for all $\alpha>0$, there exists $ B> 0 $ such that, for all $ {\mathbb {\epsilon}} $, 
    \begin{align*} 
    \left ( \int_X \left( \eta \ue^{\frac{l+1}{2}} \right)^{\crit(k)} \, dv_g  \right)^{\frac{2}{\crit(k)}} \leq \frac{K_0(n-k)(1+\alpha)}{V_m^{1-\frac{2}{\crit(k)}}} \int_X \left| {\nabla}  \right( \eta \ue^{\frac{l+1}{2}}    \left)  \right|_g^2 \, dv_g   + B \int_X \eta^2 \ue^{l+1} \, dv_g 
    \end{align*}
where $K_0(n-k)$ is as in \eqref{ineq:sobolev} and $V_m=\hbox{min}_{x\in X}\{Vol_g(Gx)/\, \hbox{dim }Gx=k\}$. By combining this inequality with \eqref{ineq:eta:2:u:alpha}, we obtain:
     \begin{equation} \label{ineq:23}
 \left(  \int_X     \left( \eta \ue^{\frac{l+1}{2}} \right)^{\critk} \, dv_g  \right)^{\frac{2}{\critk}} \chi_\eps \leq C\Vert  \ue \Vert _{l+1}^{l+1}     \end{equation}
 where $$\chi_\eps:=1 - \frac{(l+1)^2}{4l} \le       \frac{K_0(n-k)(1+\alpha)}{V_m^{1-\frac{2}{\critk}}}  \left( \int_{B_{2\delta}(Gz_0)} \ue^{\critk} \, dv_g  \right)^{1- \frac{2}{\critk}}$$
Assume that, up to extraction, 
$$\lim_{\eps\to 0}\int_{B_{2\delta}(Gz_0)} \ue^{\critk} \, dv_g<1.$$
Using \eqref{align:hyp:hme:blowup:supercritical}, there exists $1<l<\critk-1$ such that $\chi_\eps\geq \beta>0$ for all $\eps>0$ up to taking $\alpha$ small. As $\ue\to 0$ in $L^{l+1}(X)$ since $l+1<\critk$, with \eqref{ineq:23}, we then get that $\lim_{\eps\to 0}\int_{B_{\delta}(Gz_0)} \ue^{\frac{l+1}{2}\critk} \, dv_g=0$. With similar arguments, we  get that for all $\delta'<\delta$, $\ue\to 0$ in $L^q(B_{\delta'}(Gz_0)$ for all $q\geq 1$. It then follows from \eqref{principale:equation:supercritical:blowup} and elliptic theory that $\ue\to 0$ in $C^0(B_{\delta'}(Gz_0))$. Since $\int_X\ue^{\crit}\, dv_g=1$ and $X$ is compact, the existence of $x_0\in X$   such that \eqref{eq:conc} holds follows. This proves the claim.\qed

\begin{claim}
We have that $\hbox{dim }Gx_0=k$ and $\hbox{Vol}_g(Gx_0)=V_m$.
\end{claim}
\noindent{\it Proof:} We follow Faget \cite{faget2002best}. Assume that $\hbox{dim }Gx_0>k$. Therefore, there exists $\delta>0$ such that $\hbox{dim }Gx\geq k_1>k$ for all $x\in B_{2\delta}(Gx_0)$. It then follows from Hebey-Vaugon \cite{hebey0sobolev} that $H_{1,G}^2(B_\delta(Gx_0))\hookrightarrow L^p(B_\delta(Gx_0))$ is compact for $1\leq p<2^\star(k_1)$. Since $\critk<2^\star(k_1)$ and $\ue\to 0$ in $L^2(X)$, we get that $\ue\to 0$ strongly in $L^{\critk}(B_\delta(Gx_0))$, contradicting \eqref{eq:conc}. Therefore $\hbox{dim }Gx_0=k$. It follows from Faget (formula (8) in \cite{faget2002best}) that for all $\alpha>0$, there exists $ \delta_\alpha>0$ such that
\begin{equation*}
\left(\int_{B_{\delta_\alpha}(Gx_0)} |v|^{\critk}\, dv_g\right)^{\frac{2}{\critk}}\leq (1+\alpha)\frac{K_0(n-k)}{\hbox{Vol}_g(Gx_0)^{1-\frac{2}{\critk}}}\int_{B_{\delta_\alpha}(Gx_0)}|\nabla v|^2\, dv_g
\end{equation*}
for all $v\in C^1_G( B_{\delta_\alpha}(Gx_0))$ with compact support in $B_{\delta_\alpha}(Gx_0)$. Let us fix $\eta_\alpha\in C^\infty_G(X)$ with compact support in $B_{\delta_\alpha}(Gx_0)$ and such that $0\leq\eta_\alpha\leq 1$ and $\eta_\alpha(x)=1$ for $d_g(x, Gx_0)<\delta_\alpha/2$. We then get that
\begin{eqnarray}
\left(\int_{B_{\delta_\alpha/2}(Gx_0)}\ue^{\critk}\, dv_g\right)^{\frac{2}{\critk}}&\leq& \left(\int_{B_{\delta_\alpha}(Gx_0)}(\eta_\alpha\ue)^{\critk}\, dv_g\right)^{\frac{2}{\critk}}\nonumber\\
&\leq& (1+\alpha)\frac{K_0(n-k)}{\hbox{Vol}_g(Gx_0)^{1-\frac{2}{\critk}}}\int_{X}|\nabla (\eta_\alpha\ue)|^2\, dv_g\label{ineq:34}
\end{eqnarray}
Integrating by parts and using $\Vert\ue\Vert_2\to 0$, we get that
\begin{equation}
\int_{X}|\nabla (\eta_\alpha\ue)|_g^2\, dv_g=\int_X\eta_\alpha^2|\nabla\ue|_g^2\, dv_g+\int_X\eta(\Delta_g\eta)\ue^2\, dv_g\leq \int_X |\nabla\ue|_g^2\, dv_g+o(1)\label{ineq:35}
\end{equation}
Multiplying \eqref{main:eq} by $\ue$, integrating and using again $\Vert\ue\Vert_2\to 0$, we get that
\begin{equation}
\le=\int_X\le\ue^{\critk}\, dv_g=\int_X|\nabla \ue|^2\, dv_g+\int_X\he\ue^2\, dv_g=\int_X|\nabla \ue|^2\, dv_g+o(1).\label{ineq:36}
\end{equation}
Putting together \eqref{ineq:34}, \eqref{ineq:35}, \eqref{ineq:36} and \eqref{eq:conc}, we get that
$$1\leq (1+\alpha)\frac{K_0(n-k)}{\hbox{Vol}_g(Gx_0)^{1-\frac{2}{\critk}}}\le +o(1).$$
Using \eqref{align:hyp:hme:blowup:supercritical}, letting $\eps\to 0$ and then $\alpha\to 0$ yields  $\hbox{Vol}_g(Gx_0)\leq V_m$. Therefore $\hbox{Vol}_g(Gx_0)= V_m$ and  the claim is proved.\qed

     \medskip\noindent{\bf Claim: }{\it The following $L^2-$concentration holds
     \begin{equation}\label{l2:conc:g:invar}
\lim_{\eps\to 0}\frac{\int_{X \setminus B_\delta(Gx_0)   }\ue^2 \, dv_g}{\int_{X    }\ue^2 \, dv_g}  =0\hbox{ for }n-k\geq 4.
     \end{equation}
     }
We prove the claim by arguing as in Djadli-Druet \cite{dd}. We have that
 \begin{align*}
    \int_{X \setminus B_\delta(Gx_0)   }\ue^2 \, dv_g   \leq \left(\sup_{X \setminus B_\delta(Gx_0)} \ue\right)  \int_{X \setminus B_\delta(Gx_0)  } \ue  \, dv_g  
\end{align*}
Since $\ue\to 0$ in $C^0_{loc}(X\setminus Gx_0)$,  Harnack's inequality yields $ c >0$ such that  
\begin{align*}
    \int_{X \setminus B_\delta(Gx_0)  }\ue^2 \, dv_g  & \leq c \inf_{X \setminus B_\delta(Gx_0) } \ue  \int_{X \setminus B_\delta(Gx_0) }  \ue  \, dv_g  \\
    & \leq c \left( \int_{X \setminus B_\delta(Gx_0) } \ue^2 \, dv_g  \right)^{\frac{1}{2}} \int_X \ue   \, dv_g     \leq c \Vert  \ue \Vert _2 \int_X \ue   \, dv_g  
\end{align*}
Integrating \eqref{principale:equation:supercritical:blowup} yields  $ \int_X {h_{\eps }} \ue \, dv_g  = \int_X \ue^{\crit(k)-1} \, dv_g  $. It follows from \eqref{heps:intro:eng} that there exists $\beta>0$ such that $\he\geq \beta$ for all $\eps>0$. Therefore we get that
\begin{align*}
     \int_{X \setminus B_\delta(Gx_0)   }\ue^2 \, dv_g  \leq c\beta^{-1} \Vert  \ue \Vert _2  \Vert   \ue \Vert _{\crit(k)-1}^{\crit(k)-1}.
\end{align*}
if $ n-k \geq 6 $, we have  $\crit(k)-1 \leq 2 $. Using H\"older inequalities, we have   $  \Vert   \ue \Vert _{\crit(k)-1}^{\crit(k)-1} \leq c \Vert  \ue \Vert _2^{\crit(k)-1} $. Since $ \ue \to 0 $ in  $ L^2(X) $, we get that
\begin{align*}
    \int_{X \setminus B_\delta(Gx_0)  } \ue^2 \, dv_g   \leq c  \Vert  \ue \Vert _2^{\crit(k)} = o( \Vert  \ue \Vert _2^2) 
\end{align*}
if $ 4 \leq n-k \leq 5 $, we have $ \crit(k)-1 = 2 \omega + (1-\omega) \crit(k) $ with $\omega = \frac{n-k-2}{4} > 0$.  H\"older's inequality yields  $ \Vert   \ue \Vert _{\crit(k)-1}^{\crit(k)-1} \leq \Vert   \ue \Vert _{2}^{2 \omega} \Vert   \ue \Vert _{\crit(k)}^{(1-\omega)\crit(k)}$. As a result 
\begin{equation*}
     \int_{X \setminus B_\delta(Gx_0)  } \ue^2 \, dv_g   \leq   \Vert  \ue \Vert _{\crit(k)-1}^{\crit(k)-1} \Vert  \ue \Vert _2  \leq \Vert  \ue \Vert _2^{1+2\omega } \Vert  \ue \Vert _{\crit(k)}^{(1-\omega)\crit(k) } =o( \Vert  \ue \Vert _2^2) .
\end{equation*}
This proves the claim.\qed

\medskip\noindent We are now in position to take the quotient. Since $\hbox{dim }Gx_0=k$, we choose $\delta>0$ and $G'\subset G$ as in Assumption $(H)$. Then $M:=B_\delta(Gx_0)/G'$ is a manifold of dimension $n-k$ that is endowed with the metric $\bar{g}$ on $B_\delta(Gx_0)/G'$ such that the projection $(B_\delta(Gx_0),g) \to (B_\delta(Gx_0)/G',\bar{g})$ is a Riemannian submersion. We define $\bar{u}_\eps\in C^2(M)$, $\bar{h}_\eps\in C^1(M)$ and $\bar{v}\in C^2(M)$ be such that
$$\bar{u}_\eps(\bar{x})=\ue(x)\, ,\, \bar{h}_\eps(\bar{x})=\he(x)\hbox{ and }\bar{v}(\bar{x})=\hbox{Vol}_g(G'x)\hbox{ for all }x\in B_\delta(Gx_0).$$
Let us first rewrite equation \eqref{main:eq} as in Saintier \cite{saintier0blow}. Let $\bar{\varphi}\in C^\infty_c(M)$ be a function on $M=B_\delta(Gx_0)/G'$. Define $\varphi(x):=\bar{\varphi}(\bar{x})$ for all $x\in B_\delta(Gx_0)$: as one checks, $\varphi\in C^\infty_c(B_\delta(Gx_0))$ and is $G-$invariant. It then follows from \eqref{main:eq} that
$$\int_X(\nabla\ue,\nabla\varphi)_g\, dv_g+\int_X\he\ue\varphi\, dv_g=\le\int_X\ue^{\critk-1}\varphi\, dv_g.$$
We define $\tilde{g}:=\bar{v}^{\frac{2}{n-k-2}}\bar{g}$. Since $\ue,\varphi$ are $G-$invariant, we get that
$$\int_X(\nabla\ue,\nabla\varphi)_g\, dv_g=\int_{B_\delta(Gx_0)/G'}\bar{v}(\bar{x})(\nabla\bar{u}_\eps,\nabla\bar{\varphi})_{\bar{g}}\, dv_{\bar{g}}=\int_{B_\delta(Gx_0)/G'}  (\nabla\bar{u}_\eps,\nabla\bar{\varphi})_{\tilde{g}}\, dv_{\tilde{g}}.$$
Performing the same computations for the remaining terms, setting $\theps:=\bar{v}^{-\frac{2}{n-k-2}}\bar{h}_\eps$, $\tilde{f}:=\bar{v}^{-\frac{2}{n-k-2}}$ and $\tue:=\le^{\frac{1}{\crit-2}}\bar{u}_\eps$, we get that
\begin{equation*}
\Delta_{\tilde{g}}\tue+\theps\tue=  \tilde{f}\tue^{\critk-1}\hbox{ in }M .
\end{equation*}
We  deal with the $L^{\critk}-$norm. The definitions of $\tilde{f}$ and $\tilde{g}$ and \eqref{eq:conc} yield
\begin{equation*}
\lim_{\eps\to 0}\int_M \tilde{f}\tue^{\critk}\, dv_{\tilde{g}}=\lim_{\eps\to 0}\le^{\frac{\critk}{\critk-2}}=\frac{1}{K_0 (n-k)^\frac{n-k}{2}\tilde{f}(\bar{x_0})^{\frac{n-k-2}{2}}}.
\end{equation*}
Concerning the $L^2-$concentration, it follows from \eqref{l2:conc:g:invar} that for any $r<\delta$, for $n-k\geq 4$, we have that
$$\int_{B_\delta(Gx_0)\setminus B_r(Gx_0)}\ue^2\, dv_g\leq \int_{X\setminus B_r(Gx_0)}\ue^2\, dv_g=o\left(\int_{B_r(Gx_0)}\ue^2\, dv_g\right).$$
Taking the quotient, we get that
$$\int_{M\setminus B_r(\bar{x}_0)}\tue^2\, dv_{\tilde{g}} =o\left(\int_{B_r(\bar{x}_0)}\tue^2\, dv_{\tilde{g}}\right)\hbox{ when }n-k\geq 4,$$
which yields the $L^2-$concentration on $M$. We   apply Theorem \ref{theorem:Mesmar:1:article:blow:up}. Taking  $(\xe)_\eps\in X$ such that $\Vert\ue\Vert_\infty=\ue(\xe)=\me^{-\frac{n-k-2}{2}}$, we get \eqref{distance:with:orbit:smal:o:of:mu} and \eqref{ineq:co:G}. Equation \eqref{relation:de:h0:avec:courbure:scalaire:cas:variete:a:bord} rewrites
\begin{equation*}
\tilde{h}(\bar{x}_0) = \frac{n-k-2}{4(n-k-1)}\left(  \hbox{Scal}_{\tilde{g}}(\bar{x}_0)- \frac{ n-k-4}{2}\cdot \frac{ \Delta_{\tilde{g}} \tilde{f}  (\bar{x}_0)}{ \tilde{f} (\bar{x}_0) } \right)
\end{equation*}
where $\tilde{h}:=\lim_{\eps\to 0}\theps$. Using the invariance  of the conformal Laplacian, that is
$$\Delta_{\tilde{g}}\varphi+\frac{m-2}{4(m-1)} \hbox{Scal}_{\tilde{g}}\varphi=\omega^{-\frac{m+2}{m-2}}\left( \Delta_{\bar{g}}(\omega \varphi)+\frac{m-2}{4(m-1)} \hbox{Scal}_{\bar{g}}\omega \varphi\right)$$
for any $\varphi\in C^2(M)$ and where $\tilde{g}=\omega^{\frac{4}{m-2}} \bar{g}$, $m=n-k$, we get \eqref{id:scal:intro:eng}. This proves Theorem \ref{th:mesmar:2:supercitical:case:blow:up}.

\medskip\noindent{\bf Acknowledgement: } The initial version of this article required the $L^2-$concentration \eqref{hyp:conctration:L2:article:blowup} for all dimensions $n\geq 4$. The authors are grateful to the anonymous referee who noticed that this concentration could be bypassed for $n\geq 7$.

\end{document}